\documentstyle{article}
\title{Integral and Theta Formulae for Solutions of $sl_N$
Knizhnik-Zamolodchikov Equation at Level Zero}
\author{Atsushi Nakayashiki\\
Graduate School of Mathematics\\
Kyushu University}
\date{}

\begin{document}

\def\tI{\tilde{I}}
\def\bI{\bar{I}}
\def\tp{\tilde{p}}
\def\bp{\bar{p}}
\def\bk{\bar{k}}
\def\sgn{\hbox{sgn}}
\def\bcg{{\bf C}^g}
\def\bzg{{\bf Z}^g}
\def\ep{\epsilon}
\newcommand{\ch}[2]{
\left\{
\begin{array}{c}
#1\\ #2
\end{array}
\right\}_\tau}
\newcommand{\th}[2]{
\theta\left[
\begin{array}{c}
#1 \\ #2
\end{array}
\right]
}
\newcommand{\bic}[2]{
\left(
\begin{array}{c}
#1 \\ #2
\end{array}
\right)
}
\newcommand{\df}[2]{
\left(\!\!\!
\begin{array}{c}
#1\\ #2
\end{array}
\!\!\!\right)}
\newcommand{\mdf}[4]{
\left(\!\!\!
\begin{array}{c}
#1\\ #2
\end{array}
\cdots
\begin{array}{c}
#3\\ #4
\end{array}
\!\!\!\right)}

\newcommand{\amdf}[2]{
\left(\!\!\!
\begin{array}{c}
\La \\ i^1_1
\end{array}
\cdots
\begin{array}{c}
#1 \\ #2
\end{array}
\cdots
\begin{array}{c}
\La \\ i^{N-1}_{m-1}
\end{array}
\!\!\!\right)}

\newcommand{\bmdf}[4]{
\left(\!\!\!
\begin{array}{c}
#1 \\ i^1_1
\end{array}
\cdots
\begin{array}{c}
#2 \\ #3
\end{array}
\cdots
\begin{array}{c}
#4 \\ i^{N-1}_{m-1}
\end{array}
\!\!\!\right)}
\newcommand{\pdf}{
\left(\!\!\!
\begin{array}{c}
\Lambda \\ i^1_1
\end{array}
\!\!\cdots \!\Big|\!\!
\prod_{s\in I}\left(\!\!\!
\begin{array}{c}
\Lambda \\ i^r_s
\end{array}
\!\!\!\right)
\prod_{t\in J}\left(\!\!\!
\begin{array}{c}
i^{N-1}_m \\ i^{N-1}_t
\end{array}
\!\!\!\right)
\!\!\Big|\!
\cdots \!\!
\begin{array}{c}
\Lambda \\ i^{N-1}_{m-1}
\end{array}
\!\!\right)}
\newcommand{\npdf}{
\left(\!\!\!
\begin{array}{c}
\Lambda \\ i^1_1
\end{array}
\!\!\cdots \!\Big|\!\!
\prod_{s\in I}\left(\!\!\!
\begin{array}{c}
\Lambda \\ i^r_s
\end{array}
\!\!\!\right)
\prod_{t\in J}\left(\!\!\!
\begin{array}{c}
i^{N-1}_m \\ i^{N-1}_t
\end{array}
\!\!\!\right)\right)
}
\newcommand{\mmdf}[4]{
\left(\!\!\!
\begin{array}{c}
\Lambda \\ i^1_1
\end{array}
\!\!\cdots \!\Big|\!\!
\begin{array}{c}
#1 \\ i^r_1
\end{array}
\!\cdots\!\!
\begin{array}{c}
#2 \\ #3
\end{array}
\!\!\cdots\!\!
\begin{array}{c}
#4 \\ i^r_m
\end{array}
\!\!\Big|\!
\cdots \!\!
\begin{array}{c}
\Lambda \\ i^{N-1}_{m-1}
\end{array}
\!\!\right)}

\newcommand{\mmmdf}[8]{
\left(\!\!\!
\begin{array}{c}
\Lambda \\ i^1_1
\end{array}
\!\!\cdots \!\Big|\!\!
\begin{array}{c}
#1 \\ #2
\end{array}
\!\cdots\!\!
\begin{array}{c}
#3 \\ #4
\end{array}
\begin{array}{c}
#5 \\ #6
\end{array}
\!\!\cdots\!\!
\begin{array}{c}
#7 \\ #8
\end{array}
\!\!\Big|\!
\cdots \!\!
\begin{array}{c}
\Lambda \\ i^{N-1}_{m-1}
\end{array}
\!\!\right)}

\newcommand{\smdf}[4]{
\left(\!\!\!
\begin{array}{c}
\Lambda \\ i^1_1
\end{array}
\cdots\Big|
\begin{array}{c}
#1 \\ i^{N-1}_1
\end{array}
\cdots
\begin{array}{c}
#2 \\ #3
\end{array}
\cdots
\begin{array}{c}
#4 \\ i^{N-1}_{m-1}
\end{array}
\!\!\!\right)}

\newcommand{\vsmdf}[0]{
\left(\!\!\!
\begin{array}{c}
\Lambda \\ i^1_1
\end{array}
\cdots\Big|
\begin{array}{c}
i^{N-1}_m \\ i^{N-1}_1
\end{array}
\cdots
\begin{array}{c}
i^{N-1}_m \\ i^{N-1}_{m-1}
\end{array}
\!\!\!\right)}

\def\la{\lambda}
\def\La{\Lambda}
\def\kp{\kappa}
\def\ra{\longrightarrow}
\def\ot{\otimes}
\def\baf{\bar{f}}
\newtheorem{prop}{Proposition}
\newtheorem{cor}{Corollary}
\newtheorem{lem}{Lemma}
\newtheorem{theo}{Theorem}
\newtheorem{definition}{Definition}
\newtheorem{conj}{Conjecture}

\maketitle

\begin{abstract}
The solutions of the $sl_N$ Knizhnik-Zamolodchikov(KZ) equations
at level $0$ are studied. We present the integral formula
which is obtained as a quasi-classical limit of the
integral formula for form factors of the $SU(N)$ invariant
Thirring model due to F. Smirnov. 
A proof is given that those integrals satisfy $sl_N$
KZ equation of level $0$.
The relation of the integral formulae with the chiral
Szeg\"{o} kernel is clarified.
As a consequence the integral formula
with the special choice of cycles is rewritten 
in terms of the Riemann theta functions associated with 
the $Z_N$ curve. This formula gives a generalization 
of Smirnov's formula for $sl_2$.

\end{abstract}
\par

\vskip4mm
\noindent
{\Large\bf 0 \hskip4mm Introduction}
\vskip4mm
\noindent
In \cite{S2} F. Smirnov derived a curious theta formula
for the solution of the $sl_2$ Knizhnik-Zamolodchikov (KZ) 
equation at level $0$.
The aim of this paper is to generalize Smirnov's results 
to the case of $sl_N$.
Before giving a more detail of our results let us summarize
the reason why we are interested in the level $0$ case of the
KZ equation.

The KZ equation was introduced in \cite{KZ}
as one of the fundamental equations characterizing
the correlation functions of the Wess-Zumino-Witten (WZW) model
in conformal field theory.
For the affine Lie algebra $\hat{{\cal G}}$ and
its highest weight representations $V_1, \cdots, V_m$
the KZ equation has the form

\begin{eqnarray}
&&
(k+g){\partial F\over\partial \la_i}
=\sum_{j\neq i}{\Omega_{ij}\over\la_i-\la_j}F,
\nonumber
\end{eqnarray}
where $F$ is a $V_1\ot\cdots\ot V_m$ valued function in 
$\la_1,\cdots,\la_m$, $\Omega_{ij}$ is the invariant tensor,
with respect to the symmetric invariant bilinear form of ${\cal G}$,
acting on $i$-the and $j$-th tensor components, $g$ is the dual Coxeter
number of ${\cal G}$ and $k$ is a parameter.
The number $k$ is called level.
In the WZW models levels are positive integers which coinside with
those of the integrable highest weight representation of $\hat{{\cal G}}$.

The KZ equation acquires a new life from the study of
the two dimensional integrable massive quantum field theories (IMQFT)
and solvable lattice models (SLM).
F. Smirnov formulated an axiom of locality for form factors
and, for several models, obtained integral formulas of form factors \cite{S1}.
In \cite{S0} the rational $q$ deformed KZ (qKZ) equation was found as a consequence of the axiom.
Hence the moment the qKZ equation is invented the integral formula
for the solution is constructed.
It is important to note that the qKZ equation appeared in this
context is of level $0$.

Around the same time I. Frenkel and N. Reshetikhin developped
a general theory of vertex operators for quantum affine algebras \cite{FR}.
They derived a qKZ equation of general level as the equation
satisfied by the highest-highest matrix element of the vertex
operators.
This theory was successfully applied to the study of SLMs \cite{JM}.
Although, in this application to SLM, the building blocks are
the vertex operators of positive integer levels,
the form factors and the corelation functions are shown to satisfy the
level $0$ and level $-2\times(\hbox{dual Coxeter number})$ qKZ 
equation respectively \cite{JM}\cite{N2}.
Thus the level $0$ qKZ equation and their degenerations are 
of special importance in the context of IMQFT and SLM.

In order to understand the nature of form factors F. Smirnov
studied the quasi-classical limit \cite{S0}\cite{S2}.
He noticed that the period integral of the hyperelliptic curve 
$s^2=f(z)=\prod_{j=1}^{2n}(z-\la_j)$ appears as the limit
of the integral formula for the form factors of $SU(2)$ invariant
Thirring model.
Then in \cite{S2} he rewrites them in terms of Riemann theta
functions as
\begin{eqnarray}
&&
f_{\epsilon_1,\cdots,\epsilon_{2n}}(\la_1,\cdots,\la_{2n})
\nonumber
\\
&&=
\zeta_\La
(\det A)^{-3}\Delta^{-3/4}
\theta[e_\La](0)^4
\det\big(\partial_{i}\partial_{j}\log\theta[e_\La](0)\big)_{1\leq i,j\leq g},
\label{Smirtheta}
\end{eqnarray}
where $\La=(\epsilon_1,\cdots,\epsilon_{2n})$ is the sequence of $\pm$, 
the number of $+$ being equal to the number of $-$, $\zeta_\La$ a certain
fourth root of unity, 
$e_\La$ a nonsingular
even half period corresponding to the partition 
$\{1,2,\cdots,2n\}=\{j\vert \epsilon_j=+\}\sqcup \{j\vert \epsilon_j=-\}$,
$\Delta=\prod_{i<j}(\la_i-\la_j)$, 
$\partial_i=\partial/\partial z_i$, 
$g=n-1$ the genus of the curve,
$\{A_i,B_j\}$ a canonical homology basis
and $A=(\int_{A_i}z^{j-1}dz/s)_{1\leq i,j\leq n}$.
The function $F=\sum f_{\epsilon_1,\cdots,\epsilon_{2n}}
v_{\epsilon_1}\ot\cdots \ot v_{\epsilon_{2n}}$ gives a solution
to the KZ equation taking the value in $V^{\ot 2n}$ with 
$V={\bf C}v_{+}\oplus{\bf C}v_{-}$ being the vector representation of $sl_2$.

Since the theta function of an algebraic curve is the
tau function, modulo some factor, of a soliton equation,
this result suggests an intimate relation of the
level $0$ KZ equation with the soliton equations.
In spite of Smirnov's effort on this problem \cite{S2}\cite{S3}\cite{S4}
this relation is not yet clearly understood.

Integral formulas are known 
for the solutions of the KZ equation with an arbitrary level
associated with any Kac Moody Lie algebra \cite{SV1}\cite{SV2}.
In \cite{NPT} it is shown that those general integral formula has
the exact forms as their integrands in the case of $sl_2$, level $0$ and
singlet solutions.
Taking this fact into consideration is crucial to give a complete
correspondence between general formulae at level $0$ and the Smirnov type 
formulae in \cite{NPT}. A completely analogous structure exists in the
case of $sl_2$ rational qKZ equation \cite{NPT}.
Thus Smirnov type formula is related with a subtle structure of level $0$.
In the $sl_N$ case to find a similar structure to the $sl_2$ case
in the formulae in \cite{M}\cite{SV1} is not yet succeeded.

One strategy to understand Smirnov type solutions more clearly
will be to generalize it.
This is the reason why we are interested in the generalization
of the Smirnov's results to the other types of Lie algebras
than $sl_2$.

Now let us describe our results.
In \cite{S1} the integral formula for form factors of the $SU(N)$
invariant Thirring model is obtained. It is a solution
to the $sl_N$ rational qKZ equation of level zero.
We take the quasi-classical limit of this integral formula.
It is expressed as the determinant of the period integrals
of a $Z_N$ curve.
A $Z_N$ curve is a natural generalization of a hyperelliptic curve
,which corresponds to $N=2$.
Roughly speaking the integral formula obtained in this manner
should give a solution to the $sl_N$ KZ equation of level zero.
From the mathematical point of view it is not very easy to prove
rigorously that the asymptotics satisfies the KZ equation.
On the other hand the formula for the quasi-classical limit is rather simple.
Hence it is desirable and interesting 
to prove directly that it satisfies the KZ equation.
We give a proof which is new even for the $sl_2$ case. 
Compared with the proof in the generic level case 
\cite{DJMM}\cite{M}\cite{SV2} 
our proof looks more complicated. 
It will be related with the degenerate structure of Smirnov type solutions
found in \cite{NPT}.
Since we have established a correct Smirnov type
integral formula in the $sl_N$ case it is an interesting 
problem to get them from the formulae in \cite{M}\cite{SV1}
in the spirit of \cite{NPT}.

We rewrite the integral formula in terms of theta
functions on a ${\bf Z}_N$ curve.
A priori this is not a trivial task at all.
In fact the following major problems are not obvious from the
formula and arguments in the $sl_2$ case.
The first one is what kind of rational periods parametrize
the tensor components of the solution.
The second one is whether we can expect the second order derivatives
of the logarithm of theta functions or not in the $sl_N$ case.
The first problem is resolved with the help of the Thomae
formula for ${\bf Z}_N$ curves which was discovered
by Bershadsky and Radul \cite{BR}\cite{N}.
Namely the tensor component is parametrized by
certain non-singular $1/N$ or $1/2N$ periods introduced in \cite{BR}.
The second problem is solved by finding a relation
of the integrand of the integral formual with the Szeg\"{o}
kernel. In fact the product of Szeg\"{o} kernels is related
with the second order derivatives of the logarithm of theta
functions by the formula due to Fay \cite{F}.

Now the present paper is organized in the following manner.
In section $1$ the integral formula is given. 
The theta formula is given in section $2$. It is proved in section $3$.
In section 4 a proof is given that the integral formula satisfies
the KZ equation and belongs to the trivial representation of 
$sl_N$.
A derivation of fundamental relations among differential forms
used in section 4 is given in appendix.
\vskip1cm

\section{Integral Formulas}
\par
Let $sl_N$ be the simple Lie algebra of type $A_{N-1}$, $(\quad,\quad)$
the symmetric bilinear form on $sl_N$ given by $(X,Y)=\hbox{tr}(XY)$,
$\{I_j\}$ a basis of $sl_N$ and $\{I^j\}$ the dual basis with respect
to $(\quad,\quad)$. The invariant element $\Omega$ is given by
\begin{eqnarray}
&&
\Omega=\sum_j I_j\ot I^j.
\nonumber
\end{eqnarray}
Let $V$ be the $N$ dimensional irreducible representation of $sl_N$ and
$m$ a positive integer.
The Knizhnik-Zamolodchikov(KZ) equation with values in the $Nm$ fold
tensor product $V^{\ot Nm}$ of $V$ is the differential equation
for the $V^{\ot Nm}$ valued function $F$
\begin{eqnarray}
&&
(k+N){\partial F\over\partial \la_i}
=\sum_{j\neq i}{\Omega_{ij}\over\la_i-\la_j}F,
\label{KZk}
\end{eqnarray}
where $\Omega_{ij}$ means the action of $\Omega$ on the $i$-th and $j$-th
components of $V^{\ot Nm}$, $k$ is a complex number called level.
The explicit form of KZ equation in terms of the vector components
is given in section 4.

Let $v_j={}^t(0,\cdots,1,\cdots,0)$ in ${\bf C}^N$, 
where $1$ is on the $j$-th place. 
Then we have $V=\oplus_{j=1}^N{\bf C}v_j$.
We denote by $\La=(\La_1,\cdots,\La_N)$ the ordered partition
of $\{1,2,\cdots,Nm\}$ such that the number $|\La_i|$ of the elements
of $\La_i$ is $m$ for any $i$. To an ordered partition $\La$ we associate
the weight zero vector $v_\La$ of $V^{\ot Nm}$ by
\begin{eqnarray}
&&
v_{\La}=v_{k_1}\ot\cdots\ot v_{k_{Nm}},
\nonumber
\end{eqnarray}
where
\begin{eqnarray}
&&
\hbox{$i\in \La_j$ if and only if $k_i=j$}.
\nonumber
\end{eqnarray}
The set of $\{v_{\La}\}$
forms a base of the weight zero subspace of $V^{\ot Nm}$.

The operators $\sum_{j\neq i}(\la_i-\la_j)^{-1}\Omega_{ij}$
in the right hand side of (\ref{KZk}) commute with the action of $sl_N$. 
Thus it has a sense to consider the KZ equation for a function 
taking values in a fixed weight subspace of $V^{\ot Nm}$.
In this paper
we exclusively consider the solution $F$ whose value is
in the weight zero subspace of the tensor product $V^{\ot Nm}$.
Then we can define the component $f_\La$ of $F$ by
\begin{eqnarray}
&&
F=\sum_{\La}f_{\La}v_{\La},
\label{component}
\end{eqnarray}
where the sum is over all ordered partition $\La$. 

We denote by $C$ the compact Riemann surface defined from the equation 
$s^N=f(z)=\prod_{j=1}^{Nm}(z-\la_i)$.
It is called a $Z_N$ curve \cite{BR,N}.
The genus $g$ of $C$ is given by $2g=(N-1)(Nm-2)$.
For $\La_r$ and $p\in\La_r$ set
\begin{eqnarray}
&&
g^{(\La_r)}(z)=\prod_{j\not\in\La_r}(z-\la_j),\ \
g_{\La_r}(z)=\prod_{j\in\La_r}(z-\la_j),\ \
g^{(p)}_{\La_r}(z)=\prod_{j\in\La_r,j\neq p}(z-\la_j)
\nonumber
\end{eqnarray}
and define the meromorphic differential form $\mu_p^\La(z)$ on $C$ by
\begin{eqnarray}
&&
\mu_p^\La(z)=
{g^{(\La_r)}(\la_p)g^{(p)}_{\La_r}(z) \over (z-\la_p)s }dz.
\nonumber
\end{eqnarray}
We set $L=(N-1)m-1$.
Then we have

\begin{theo}\label{intformula}
Let $\{p_1,\cdots,p_L\}$ is an arbitrary subset of
$\{1,2,\cdots,Nm\}$. Define
\begin{eqnarray}
&&
f(\la_1,\cdots,\la_{Nm})_\La
=
{\Delta^{{N-1\over N^2}}
\over
\prod_{i<j}(\La_i\La_j)}
{\det(\int_{\gamma_i}\mu_{p_j}^\La)_{1\leq i,j\leq L},
\over
\Delta(p_1,\cdots,p_L)}
\label{intsol}
\end{eqnarray}
where
$(\La_i\La_j)=\prod_{r\in\La_i,s\in\La_j}(\la_r-\la_s)$,
$\Delta(p_1,\cdots,p_L)=\det(\la_{p_j}^{L-i})_{1\leq i,j\leq L}$ 
and $\Delta=\prod_{i<j}(\la_i-\la_j)$.
Then
\begin{description}
\item[(0)]
The right hand side of (\ref{intsol}) does not depend on the choice of
$\{p_1,\cdots,p_L\}$.
\item[(1)]
The function $F$ given by (\ref{component}) and (\ref{intsol})
is a solution to the $sl_N$ KZ equation of level zero 
for arbitrary set of $L$ cycles $\{\gamma_1,\cdots,\gamma_L\}$ 
on $C$.
\item[(2)]
For any $X\in sl_N$, $XF=0$.
\end{description}
\end{theo}
\vskip5mm

The first statement of Theorem \ref{intformula} follows from another expression
for $f_\La$.
Let us set
\begin{eqnarray}
&&
\zeta_j^{\La}={dz\over s}\sum_{k=1}^Ng_{\La_k}(z)
\Big[
{d\over dz}
{g^{(\La_k)}(z)\over z^{L-j+1}}
\Big]_0,
\nonumber
\end{eqnarray}
where $[\quad]_0$ denotes the polynomial part of a Laurent polynomial.
It is obvious that $d/dz$ can be out side of the symbol $[\quad]_0$.

\begin{theo}\label{intformula2}
The function $f_\La$ given by (\ref{intsol}) is also written as
\begin{eqnarray}
&&
f(\la_1,\cdots,\la_{Nm})_\La
=
{\Delta^{{N-1\over N^2}}
\over
\prod_{i<j}(\La_i\La_j)}
\det(\int_{\gamma_i}\zeta^{\La}_j)_{1\leq i,j\leq L}.
\label{intsol2}
\end{eqnarray}
\end{theo}
\vskip8mm

We shall give some comments on the integral formula given here.
In \cite{S1} F. Smirnov derived the integral formula of
form factors of $SU(N)$ invariant Thirring model which satisfy
the deformed Knizhnik-Zamolodchikov(dKZ) equation on level zero. 
Scaling the rapidity variables $\beta_j$ as $\beta_j=\la_j/h$ and 
taking the quasi-classical limit $h\rightarrow 0$, 
we obtain the integral formula in Theorem \ref{intformula2}
with some special choice of cycles $\{\gamma_i\}$.

\section{Theta Formula}
\par
We shall give another expression for the solution $F$ given in
Theorem \ref{intformula}. To give a precise statement we
prepare necessary notations associated with the $Z_N$ curve $C$ \cite{N}.
The $N$-cyclic automorphism $\phi$ of $C$ is defined by
$\phi:(z,s)\mapsto(z,\omega s)$, where $\omega$ is the $N$-th
primitive root of unity.
There are $Nm$ branch points $Q_1,\cdots,Q_{Nm}$ whose projection
to $z$ coordinate are $\la_1,\cdots,\la_{Nm}$.
The basis of holomorphic 1-forms on $C$ is given by
\begin{eqnarray}
&&
w^{(\alpha)}_\beta={z^{\beta-1}dz\over s^{\alpha}}
\quad
1\leq\alpha\leq N-1,
\quad
1\leq\beta\leq \alpha m-1.
\nonumber
\end{eqnarray}

We fix a canonical homology basis $\{\alpha_i,\beta_j\}$
whose intesection numbers are 
$\alpha_i\cdot\alpha_j=\beta_i\cdot\beta_j=0$, 
$\alpha_i\cdot\beta_j=\delta_{ij}$.
Let $\Delta$ be a Riemann divisor 
for this choice of canonical basis. Let us define the divisor class
$D$ by $D=NQ_i$ which is independent of the choice of $i$.
To each ordered partition $\La$ we associate the divisor class
$e_\La$ \cite{N} by
\begin{eqnarray}
&&
e_\La \equiv \La_1+2\La_2+\cdots+(N-1)\La_{N-1}-D-\Delta,
\nonumber
\end{eqnarray}
where for a subset $S$ of $\{1,2,\cdots, Nm\}$
we set
\begin{eqnarray}
&&
S=\sum_{j\in S}Q_j.
\nonumber
\end{eqnarray} 
The divisor class $e_\La$ is a $1/N$ period for $N$ even and is a
$1/2N$ period for $N$ odd.
We consider the index of $\La_j$ by modulo $N$. In particular
$\La_0=\La_N$.

Let $\{v_j(x)\}$ be the basis of the normalized abelian differentials
of the first kind whose normalization is
\begin{eqnarray}
&&
\int_{A_j}v_k(x)=2\pi i\delta_{jk}.
\nonumber
\end{eqnarray}
We set $\displaystyle{\tau_{jk}=\int_{B_j}v_k(x)}$.
Then the period matrix $\tau=(\tau_{jk})$ is symmetric and its real part
is negative definite. The Jacobian variety $J(C)$ of $C$ is
described as $J(C)={\bf C}^g/2\pi i{\bf Z}^g+{\bf Z}^g\tau$.
For any elememt $e\in{\bf C}^g$, there exist unique elememts 
$\delta,\epsilon\in{\bf R}^g$ such that
\begin{eqnarray}
&&
e=\ch{\delta}{\epsilon}=2\pi i\epsilon+\delta\tau.
\nonumber
\end{eqnarray}
We call $\delta,\epsilon$ the characteristics of $e$.
The Riemann theta function with characteristics $\delta,\epsilon$
is defined by
\begin{eqnarray}
&&
\th{\delta}{\epsilon}(z)=
\sum_{m\in{\bf Z}^g}\exp
\Big(
{1\over2}(m+\delta)\tau(m+\delta)^t+(z+2\pi i\epsilon)(m+\delta)^t
\Big).
\nonumber
\end{eqnarray}
It satisfies the equation
\begin{eqnarray}
\th{\delta+m}{\epsilon+n}(z)=
\exp(2\pi in\delta^t)\th{\delta}{\epsilon}(z),
\label{charshift}
\end{eqnarray}
for $m,n\in{\bf Z}^g$.
For an ordered partition $\La$ let us take a representative 
$\bar{e}_\La\in{\bf C}^g$ of $e_\La$ and let 
$\displaystyle{\bar{e}_\La=\ch{\delta}{\epsilon}}$. 
Then the logarithmic derivatives
\begin{eqnarray}
&&
\partial^\alpha\log\th{\delta}{\epsilon}(z),
\quad |\alpha|\geq1,
\nonumber
\end{eqnarray}
are independent of the choice of the representative $\bar{e}_\La$
by (\ref{charshift}),
where $\alpha=(\alpha_1,\cdots,\alpha_g)$, 
$|\alpha|=\alpha_1+\cdots+\alpha_g$,
$\partial^\alpha
=\partial_{1}^{\alpha_1}\cdots\partial_{g}^{\alpha_g}$
and $\partial_j=\partial/\partial z_j$.
Hence we use the notation $\partial^\alpha\log\theta[e_\La](z)$
for those logarithmic derivatives for the sake of simplicity.

Let us define the connection matrix between 
$\{w^{(\alpha)}_\beta\}$ and $\{v_j(x)\}$ by
\begin{eqnarray}
&&
v_j(x)=\sum_{\alpha,\beta}\sigma_{j(\alpha\beta)}w^{(\alpha)}_\beta.
\nonumber
\end{eqnarray}
With the aid of this matrix we define the vector
field on $J(C)$ by
\begin{eqnarray}
&&
D_\beta=
\sum_{j=1}^g\sigma_{j(N-1\beta)}
\partial_j,
\quad
1\leq \beta\leq L.
\nonumber
\end{eqnarray}

Now we can state the theta formula.

\begin{theo}\label{theta1}
For any subset $\{i_1<\cdots<i_L\}$ of $\{1,2,\cdots,g\}$
we take the cycles $\{\gamma_j\}$ as $\gamma_j=A_{i_j}$.
Then the corresponding solution of the KZ equation 
in Theorem \ref{intformula} is given by
\begin{eqnarray}
&&
f(\la_1,\cdots,\la_{Nm})_\La=
c
{\Delta^{{N-1\over N^2}}
\over
\prod_{i<j}(\La_i\La_j)}
\det\big(\partial_{i_j}D_k\log\theta[e_\La](0)\big)_{1\leq j,k\leq L},
\nonumber
\end{eqnarray}
where $c$ is the overall constant independent of $\la_i$'s and $\La$.
\end{theo}
\vskip5mm

Using the Thomae formula for ${\bf Z}_N$ curves one can rewrite
the $\prod_{i<j}(\La_i\La_j)$ in terms of theta constants.
The result is

\begin{theo}\label{theta2}
For the same choice of cycles as in Theorem \ref{theta1},
we have
\begin{eqnarray}
&&
f(\la_1,\cdots,\la_{Nm})_\La=
C(\la)\zeta_\La
\prod_{\sigma\in S_{N-1}}
\theta[e_{\La^\sigma}](0)^{{12N\over(N+1)!}}
\det\big(\partial_{i_j}D_k\log\theta[e_\La](0)\big)_{1\leq j,k\leq L},
\nonumber
\end{eqnarray}
where $\zeta_\La$ is some $N(N+1)!/3$-th root of unity,
$C(\la)$, which is independent of the partition $\La$,
is given by
\begin{eqnarray}
&&
C(\la)=c(\det A)^{-{6\over N+1}}
\Delta^{-3{N-1\over N+1}+{N-1\over N^2}}.
\nonumber
\end{eqnarray}
Here $c$ is a constant independent of $\la_i$'s
and  $\displaystyle{A=\det(\int_{A_i}w^{(\alpha)}_\beta)}$.
For an elememt $\sigma$ of the symmetric group $S_{N-1}$ of
degree $N-1$ we define
\begin{eqnarray}
&&
\La^\sigma=(\La_0,\La_{\sigma(1)},\cdots,\La_{\sigma(N-1)}).
\nonumber
\end{eqnarray}
\end{theo}
\vskip5mm

\noindent
{\bf Remark.} If $N=2$, then $L=m-1=g$ and $(\sigma_{j(1k)})=A^{-1}$,
where $g$ is the genus of the hyperelliptic curve $C$.
In particular $(i_1,\cdots,i_L)=(1,\cdots,g)$.
From the matrix relation
\begin{eqnarray}
&&
\big(\partial_{i}D_k\log\theta[e_\La](0)\big)=
\big(\partial_{i}\partial_{j}\log\theta[e_\La](0)\big)A^{-1}
\nonumber
\end{eqnarray}
we have
\begin{eqnarray}
&&
f_\La=c'\zeta_\La
(\det A)^{-3}\Delta^{-3/4}
\theta[e_\La](0)^4
\det\big(\partial_{i}\partial_{j}\log\theta[e_\La](0)\big)_{1\leq i,j\leq g}
\nonumber
\end{eqnarray}
which is nothing but the Smirnov's formula (\ref{Smirtheta}) for $sl_2$.
\vskip5mm


\section{Proof of the Theta Formula}
Let $\tilde{C}$ be the universal covering space of $C$.
We identify a holomorphic one forms on $C$ with those on 
$\tilde{C}$ which are invariant under the action of the
fundamental group of $C$. We set $v=(v_1,\cdots,v_g)$, the
vector of the normalized differentials of the first kind.
Recall that the chiral Szeg\"{o} kernel defined by $e_\La$
is
\begin{eqnarray}
&&
R(x,y\vert e_\La)=
{
\theta[e_\La](y-x)\over\theta[e_\La](0)E(x,y)
}
\qquad x,y\in \tilde{C},
\nonumber
\end{eqnarray}
where 
\begin{eqnarray}
&&
y-x=\int_{x}^{y}v,
\nonumber
\end{eqnarray}
the integral being taken in $\tilde{C}$ and $E(x,y)$ 
is the prime form \cite{N}.
We remark that, as to the $e_\La$ dependence, $R(x,y\vert e_\La)$ depends
only on the divisor class of $e_\La$.

Let $\omega(x,y)$ be the canonical symmetric differential, 
that is, $\omega(x,y)$ is the section of the canonical bundle
of $C\times C$, symmetric in $x$ and $y$, has a double pole
at $x=y$, has the vanishing $A$ period in each of the variables
and has some normalization (for more precise definition see \cite{N}).
Theorem \ref{theta1} is a corollary of the following proposition.

\begin{prop}\label{muptheta}
For $1\leq p\leq Nm$ we have
\begin{eqnarray}
&&
\mu_p^\La=Nf'(\la_p)^{{N-1\over N}}\omega(x,Q_p)+
N^2\sum_{i=1}^g\sum_{\beta=1}^L
\la_p^{\beta-1}D_\beta\partial_i\log\theta[e_\La](0)v_i(x).
\nonumber
\end{eqnarray}
\end{prop}
\vskip5mm

As in \cite{N} the value of a (half)differential form
at the branch point $Q_p$ is defined as the coefficient
of $dt$ (or $\sqrt{dt}$ in the half differential case)
in the expansion of the form in the local coordinate
$t=(z-\la_p)^{1/N}$.
Assuming this proposition let us first prove Theorem \ref{theta1}.
\vskip3mm

\noindent
Proof of Theorem \ref{theta1}:
\par
\noindent
Since the integral of $\omega(x,Q_p)$ along the cycle $A_i$ is zero
for any $i$, we have
\begin{eqnarray}
&&
\int_{A_i}\mu_p^\La=
N^2\sum_{\beta=1}^L\la_p^{\beta-1}D_\beta\partial_i\log\theta[e_\La](0),
\nonumber
\end{eqnarray}
where we use the normalization condition of $\{v_j\}$.
Thus we have
\begin{eqnarray}
&&
\det\Big(\int_{A_{i_j}}\mu_{p_k}^\La\Big)_{1\leq j,k\leq L}
\nonumber
\\
&=&
N^{2L}\det
\big(\partial_{i_j}D_\beta\log\theta[e_\La](0)\big)_{1\leq j,\beta\leq L}
\det\big(\la_{p_l}^{k-1})_{1\leq k,l\leq L}
\nonumber
\\
&=&
(-1)^{{L(L-1)\over2}}N^{2L}\det
\big(\partial_{i_j}D_\beta\log\theta[e_\La](0)\big)_{1\leq j,\beta\leq L}
\det\big(\la_{p_l}^{L-k})_{1\leq k,l\leq L}
\nonumber
\\
&=&
(-1)^{{L(L-1)\over2}}N^{2L}\det
\big(\partial_{i_j}D_\beta\log\theta[e_\La](0)\big)_{1\leq j,\beta\leq L}
\Delta(p_1,\cdots,p_L).
\nonumber
\end{eqnarray}
Substituting this equation into (\ref{intsol}) we obtain the formula
of Theorem \ref{theta1}. $\Box$
\vskip3mm

In order to prove Proposition \ref{muptheta}
we first prove

\vskip5mm
\begin{prop}\label{key}
For $1\leq p\leq Nm$ we have
\begin{eqnarray}
&&
\mu_p^\La=Nf'(\la_p)^{{N-1\over N}}
R(x,Q_p\vert e_\La)R(x,Q_p\vert -e_\La).
\label{doubleR}
\end{eqnarray}
\end{prop}
\vskip5mm

For the proof of this proposition 
let us recall the following notation\cite{N,BR}:
\begin{eqnarray}
{\cal L}&=&\{-{N-1\over2},-{N-1\over2}+1,\cdots,{N-1\over2}\},
\nonumber
\\
q_l(i)&=&{1-N\over 2N}+\Big\{{l+i+{N-1\over2}\over N}\Big\},
\nonumber
\end{eqnarray}
where  $l\in {\cal L}+{\bf Z}$, $i\in {\bf Z}$ and  $\{a\}=a-[a]$ 
is the fractional part of 
$a\in {\bf Q}$.

In \cite{N} we have proved

\begin{prop}
For an ordered partition $\La$ we have
\begin{eqnarray}
R(x,y\vert e_\La)&=&{1\over N}
{\sum_{l\in{\cal L}}f_l(x,\La)f_{-l}(y,\La^{-})
\over z(y)-z(x)}
\label{szegoexp}
\\
f_{\pm l}(x,\La^{\pm})&=&
\prod_{i=1}^{Nm}(z(x)-\la_i)^{\pm q_l(k_i)}\sqrt{dz(x)},
\nonumber
\end{eqnarray}
where $\La^{-}=(\La_0,\La_{N-1},\cdots,\La_1)$, $\La^{+}=\La$
and $k_i=j$ is determined by $i\in\La_j$.
\end{prop}
\vskip5mm

By Proposition 4 in \cite{N}, for $l=-(N-1)/2+j$ with
$j\in{\bf Z}$, we have
\begin{eqnarray}
\hbox{div}f_l(x,\La^{\pm})&=&\La_{\pm 1\mp j}+2\La_{\pm 2\mp j}
+\cdots+(N-1)\La_{\mp 1\mp j}
-\sum_{k=1}^N\infty^{(k)}.
\label{div}
\end{eqnarray}
where $\{\infty^{(k)}\}$ are $N$ infinity points.
The index $k$ of $\La_k$ is considered by modulo $N$.
The $\La_{\mp j}$ is missing in the right hand side of (\ref{div}).
Let $p\in\La_r$, $0\leq r\leq N-1$.
We set $y=Q_p$ in (\ref{szegoexp}). 
By (\ref{div}) only the term 
$-l=-(N-1)/2+r$ in the sum of the
right hand side of (\ref{szegoexp}) is alive.
Hence we have
\begin{eqnarray}
&&
R(x,Q_p\vert e_\La)={1\over N}
{f_{-{N-1\over2}+r^{-}}(x,\La)f_{-{N-1\over2}+r}(Q_p,\La^{-})
\over \la_p-z(x)},
\nonumber
\end{eqnarray}
where $r^{-}=N-1-r$.
Similarly 
\begin{eqnarray}
&&
R(x,Q_p\vert -e_\La)=
R(x,Q_p\vert e_{\La^{-}})=
{1\over N}
{f_{-{N-1\over2}+r-1}(x,\La^{-})f_{-{N-1\over2}+N-r}(Q_p,\La)
\over \la_p-z(x)}.
\nonumber
\end{eqnarray}
Therefore
\begin{eqnarray}
&&
R(x,Q_p\vert e_\La)R(x,Q_p\vert -e_\La)=
\nonumber
\\
&&
{1\over N^2\big(z(x)-\la_p\big)^2}
\Big(
f_{-{N-1\over2}+r^{-}}(x,\La)
f_{-{N-1\over2}+r-1}(x,\La^{-})
f_{-{N-1\over2}+r}(Q_p,\La^{-})
f_{-{N-1\over2}+N-r}(Q_p,\La)
\Big)
\label{prodrr}
\end{eqnarray}
Let us calculate the right hand side of (\ref{prodrr}).

\begin{lem}
The following equations hold:
\begin{eqnarray}
f_{-{N-1\over2}+r}(Q_p,\La^{-})
f_{-{N-1\over2}+N-r}(Q_p,\La)
&=&
{Nf'(\la_p)^{{1\over N}}\over g^{(p)}_{\La_r}(\la_p)},
\label{prod1}
\\
f_{-{N-1\over2}+r^{-}}(x,\La)
f_{-{N-1\over2}+r-1}(x,\La^{-})
&=&
{g_{\La_r}(z)\over s}dz.
\label{prod2}
\end{eqnarray}
\end{lem}
\vskip3mm

\noindent
Proof. Let $t=(z-\la_p)^{1/N}$ be the local coordinate around $Q_p$.
Then
\begin{eqnarray}
f_{-{N-1\over2}+r}(x,\La^{-})
&=&
\sqrt{N}\prod_{j\neq p}^{Nm}
(\la_p-\la_j)^{q_{-{N-1\over2}+r}(N-k_j)}
\sqrt{dt}(1+O(t^N)),
\nonumber
\\
f_{-{N-1\over2}+N-r}(x,\La)
&=&
\sqrt{N}\prod_{j\neq p}^{Nm}
(\la_p-\la_j)^{q_{-{N-1\over2}+N-r}(k_j)}
\sqrt{dt}(1+O(t^N)),
\nonumber
\end{eqnarray}
where we use the relation $-q_l(i)=q_{-l}(N-i)$ \cite{N}.
Therefore we need to calculate the number
\begin{eqnarray}
&&
q_{-{N-1\over2}+r}(N-k_j)+q_{-{N-1\over2}+N-r}(k_j).
\nonumber
\end{eqnarray}
By a direct calculation we have
\begin{eqnarray}
q_{-{N-1\over2}+r}(N-i)+q_{-{N-1\over2}+N-r}(i)=
\left\{
\begin{array}{ll}
-1+{1\over N} & \mbox{if $i=r$ mod $N$}\\
{1\over N}   & \mbox{if $i\neq r$ mod $N$.}
\end{array}
\right.
\nonumber
\end{eqnarray}
Hence
\begin{eqnarray}
&&
f_{-{N-1\over2}+r}(Q_p,\La^{-})
f_{-{N-1\over2}+N-r}(Q_p,\La)
\nonumber
\\
&=&
N\prod_{j\notin\La_r}(\la_p-\la_j)^{1/N}
\prod_{j\in\La_r\backslash\{p\}}(\la_p-\la_j)^{-1+1/N}
\nonumber
\\
&=&
{Nf'(\la_p)^{1/N} \over
\prod_{j\in\La_r\backslash\{p\}}(\la_p-\la_j)}.
\nonumber
\end{eqnarray}
Similarly, using
\begin{eqnarray}
&&
q_{-{N-1\over2}+r^{-}}(i)+q_{-{N-1\over2}+r-1}(N-i)
=
\left\{
\begin{array}{ll}
1-{1\over N} & \mbox{if $i=r$ mod $N$}\\
-{1\over N}   & \mbox{if $i\neq r$ mod $N$}
\end{array}
\right.
\nonumber
\end{eqnarray}
we have
\begin{eqnarray}
&&
f_{-{N-1\over2}+r^{-}}(x,\La)
f_{-{N-1\over2}+r-1}(x,\La^{-})
=
{\prod_{j\in\La_r}(z-\la_j)\over s}dz.
\nonumber
\end{eqnarray}
Thus the lemma is proved. $\Box$
\vskip3mm

Multiplying (\ref{prod1}) and (\ref{prod2}) 
and substituting it into (\ref{prodrr})
we obtain the equation (\ref{doubleR}). 
$\Box$
\vskip8mm

Recall the Fay's formula (\cite{F}, Corollary 2.12, see also \cite{N}
section 4) :
\begin{eqnarray}
&&
R(x,y\vert e_\La)R(x,y\vert -e_\La)
=
\omega(x,y)+
\sum_{i,j=1}^g
{\partial^2\log\theta[e_\La]\over\partial z_i\partial z_j}
(0)v_i(x)v_j(y).
\label{remarkable}
\end{eqnarray}
By calculation we have
\begin{eqnarray}
&&
v_j(Q_p)={N\over f'(\la_p)^{1-{1\over N}}}
\sum_{\beta=1}^L\sigma_{j(N-1\beta)}\la_p^{\beta-1}.
\label{vjexp}
\end{eqnarray}
Substituting (\ref{vjexp}) into
(\ref{remarkable}) and using Proposition \ref{key} 
we have the equation in Proposition \ref{muptheta}. $\Box$
\vskip5mm

In order to prove Theorem \ref{theta2} let us recall the Thomae
formula for ${\bf Z}_N$ curves \cite{BR}\cite{N}:
\begin{eqnarray}
&&
\theta[e_\La](0)^{2N}
=C_{\La}(\det A)^{N}\prod_{i\leq j}(\La_i\La_j)^{2Nq(i,j)+N\mu},
\nonumber
\end{eqnarray}
where $C_\La$ is a constant independent of $\la_i$'s and
\begin{eqnarray}
&&
q(i,j)=\sum_{l\in{\cal L}}q_l(i)q_l(j),
\quad
\mu={(N-1)(2N-1)\over 6N},
\nonumber
\\
&&
(\La_i\La_j)=\prod_{r\in\La_i,s\in\La_j}(\la_r-\la_s)
\quad\hbox{for $i\neq j$},
\nonumber
\\
&&
(\La_i\La_i)=\prod_{r,s\in\La_i,r<s}(\la_r-\la_s).
\end{eqnarray}
The number $q(i,j)$ depends only on $\vert i-j\vert$ \cite{N}.
In particular $q(0,0)=q(i,i)$, $q(i,j)=q(j,i)$
for any $i$ and $j$.
We define the action of the symmetric group $S_{N}$ of degree $N$
on the set of ordered partitions by
\begin{eqnarray}
&&
\La^\sigma=(\La_{\sigma(0)},\cdots,\La_{\sigma(N-1)}),
\quad \sigma\in S_N.
\nonumber
\end{eqnarray}
The subgroup $S_{N-1}$ acts on the index $1,2,\cdots,N-1$ as we
already defined.

\begin{prop}\label{Thomae1}
For an ordered partition $\La$ we have
\begin{eqnarray}
&&
\prod_{i<j}(\La_i\La_j)=
C\bar{\zeta}_\La(\det A)^{{6\over N+1}}
\Delta^{3{N-1\over N+1}}
\theta[e_{\La^\sigma}](0)^{-{12N\over(N+1)!}},
\nonumber
\end{eqnarray}
where $\bar{\zeta}_\La$ is some ${N(N+1)!/3}$-th root of unity
and $C$ is a constant independent of $\la_i$'s and $\La$.
\end{prop}
\vskip5mm

Let us prove Proposition \ref{Thomae1}.
By taking the product of theta function with the characteristics
$e_{\La^\sigma}$ for all $\sigma\in S_{N}$ we have
\begin{eqnarray}
\prod_{\sigma\in S_{N}}\theta[e_{\La^\sigma}](0)^{2N}
&=&
\big(\prod_{\sigma\in S_{N}}C_{\La^\sigma}\big)
(\det A)^{N!N}
\prod_{\sigma\in S_{N}}\prod_{i\leq j}
(\La_{\sigma(i)}\La_{\sigma(j)})^{2Nq(i,j)+N\mu}
\nonumber
\\
&=&
\pm
\big(\prod_{\sigma\in S_{N}}C_{\La^\sigma}\big)
(\det A)^{N!N} \Delta^{N!N\mu}
\prod_{i\leq j}
(\La_{i}\La_{j})^{2N
\sum_{\sigma\in S_{N}}q(\sigma(i),\sigma(j))}.
\nonumber
\end{eqnarray}
Set $\gamma=q(0,0)$ and
\begin{eqnarray}
F&:=&\prod_{i\leq j}(\La_{i}\La_{j})^{2N
\sum_{\sigma\in S_{N}}q(\sigma(i),\sigma(j))}
\nonumber
\\
&=&
\prod_{i=0}^{N-1}(\La_i\La_i)^{N!\cdot 2N\gamma}
\prod_{i<j}(\La_{i}\La_{j})^{2N
\sum_{\sigma\in S_{N}}q(\sigma(i),\sigma(j))}.
\nonumber
\end{eqnarray}
Since, for $i\neq j$,
\begin{eqnarray}
&&
\sum_{\sigma\in S_{N}}q(\sigma(i),\sigma(j))=
2\cdot(N-2)!\sum_{r<s}q(r,s),
\nonumber
\end{eqnarray}
we have
\begin{eqnarray}
&&
F=\prod_{i=0}^{N-1}(\La_i\La_j)^{N!\cdot 2N\gamma}
\prod_{i<j}(\La_{i}\La_{j})^{4N\cdot(N-2)!
\sum_{i<j}q(i,j)}.
\nonumber
\end{eqnarray}
Using $q(i,j)=q(i+1,j+1)$ and Lemma 10 in \cite{N}
we have
\begin{eqnarray}
&&
\sum_{i<j}q(i,j)=-{N^2-1\over24},
\quad
\gamma={N^2-1\over12N}.
\nonumber
\end{eqnarray}
Thus 
\begin{eqnarray}
&&
F=\Delta^{{(N+1)!(N-1)\over6}}
\prod_{i<j}(\La_{i}\La_{j})^{-{(N+1)!N\over6}}.
\nonumber
\end{eqnarray}
From this we obtain
\begin{eqnarray}
&&
\prod_{\sigma\in S_{N}}\theta[e_{\La^\sigma}](0)^{2N}
=\pm
\big(\prod_{\sigma\in S_{N}}C_{\La^\sigma}\big)
(\det A)^{N!N} \Delta^{{N!N(N-1)\over2}}
\prod_{i<j}(\La_{i}\La_{j})^{-{(N+1)!N\over6}}.
\nonumber
\end{eqnarray}
Recall that the ordered partitions which are obtained
from $\La$ by the cyclic permutation of indices correspond
to linear equivalent divisor $e_\La$ \cite{N}.
Therefore 
\begin{eqnarray}
&&
\prod_{\sigma\in S_{N-1}}\theta[e_{\La^\sigma}](0)^{2N^2}
=\pm
\big(\prod_{\sigma\in S_{N-1}}C_{\La^\sigma}^N\big)
(\det A)^{N!N} \Delta^{{N!N(N-1)\over2}}
\prod_{i<j}(\La_{i}\La_{j})^{-{(N+1)!N\over6}}.
\nonumber
\end{eqnarray}
Since $C_{\La}^{2N}$ does not depend on $C_\La$,
we have the equation in Proposition \ref{Thomae1}. $\Box$
\vskip1cm

In this section we shall give a proof of Theorem \ref{intformula} and 
Theorem \ref{intformula2}.

\subsection{Proof of Theorem 2}  
Theorem \ref{intformula2} follows from the following proposition.

\begin{prop}\label{exact}

For any $p$
\begin{eqnarray}
&&
\sum_{j=1}^L\zeta_j^{\La}(z)\la_p^{L-j}
=\mu^\La_p(z)+Nd{s^{N-1}\over z-\la_p}.
\nonumber
\end{eqnarray}

\end{prop}

The proof of this proposition is totally
similar to the case of $sl_2$ \cite{S2}.
For the sake of making the paper selfcontained we give a proof. 

\noindent
{\it Proof.} Let $p\in \La_r$. It is sufficient to prove
the following equation, the coefficient of $dz/s$ :
\begin{eqnarray}
&&
\sum_{j=1}^L\sum_{k=1}^N g_{\La_k}(z)
\Big[{d\over dz}
{g^{(\La_k)}(z)\over z^{L-j+1}}
\Big]_0\la_p^{L-j}
\nonumber
\\
&&
={g^{(\La_r)}(\la_p)g^{(p)}_{(\La_r)}(z)\over z-\la_p}
-{Nf(z)\over (z-\la_p)^2}+{(N-1)f'(z)\over z-\la_p}.
\label{1eq}
\end{eqnarray}
Since both hand sides are rational functions in $z$, it is sufficient to prove
(\ref{1eq}) for $|z|$ sufficiently large.

Let $t$ be a complex parameter.
More generally we calculate the right hand side of (\ref{1eq}) 
replaced $\la_p$ by $t$.
We have
\begin{eqnarray}
\sum_{j=1}^L\Big[{d\over dz}
{g^{(\La_k)}(z)\over z^{L-j+1}}
\Big]_0t^{L-j}
&=&
{d\over dz}\Big[\sum_{j=0}^\infty(tz^{-1})^jz^{-1}g^{(\La_k)}(z)
\Big]_0
\nonumber
\\
&=&
{d\over dz}\Big[{g^{(\La_k)}(z)\over z-t}\Big]_0
\nonumber
\\
&=&
{d\over dz}{g^{(\La_k)}(z)-g^{(\La_k)}(t)\over z-t}.
\nonumber
\end{eqnarray}
Thus 
\begin{eqnarray}
&&\sum_{j=1}^L\sum_{k=1}^N
g_{\La_k}(z)
\Big[{d\over dz}{g^{(\La_k)}(z)\over z^{L-j+1}}\Big]_0t^{L-j}
\nonumber
\\
&&=
-{Nf(z)\over (z-t)^2}+{(N-1)f'(z)\over z-t}
+\sum_{k=1}^N{g_{\La_k}(z)g^{(\La_k)}(t)\over (z-t)^2}.
\label{2eq}
\end{eqnarray}
In this calculation we use
\begin{eqnarray}
&&
\sum_{k=1}^Ng_{\La_k}(z){d\over dz}g^{(\La_k)}(z)=(N-1)f'(z).
\nonumber
\end{eqnarray}
If we set $t=\la_p$, $p\in \La_r$ in (\ref{2eq}), then we get 
the right hand side of (\ref{1eq}). $\Box$
\vskip5mm

\subsection{Some Notations}

For an ordered partition $\La$ let us denote by $\La^{(ij)}$
the ordered partition which is obtained from $\La$ by
exchanging $i$ and $j$. 
For example if $\La_i=\{(i-1)m+1,\cdots,im\}$, $1\leq i\leq N$,
then $\La^{(12)}=\La$ (assuming $m\geq 2$), 
$\La^{(1m+1)}_1=\{m+1,2,\cdots,m\}$,
$\La^{(1m+1)}_2=\{1,m+2,\cdots,2m\}$, 
$\La^{(1m+1)}_k=\La_k$ for $k\geq 3$.

In terms of the components $f_\La$, if $p\in \La_i$,
the KZ equation of level zero is

\begin{eqnarray}
&&
N{\partial f_\La\over\partial \la_p}
=\Big(
(1-{1\over N})\sum_{j\in \La_i,j\neq p}{1\over\la_p-\la_j}
-{1\over N}\sum_{j\not\in \La_i}{1\over\la_p-\la_j}
\Big)f_\La
+\sum_{j\not\in \La_i}{1\over\la_p-\la_j}f_{\La^{(pj)}}.
\nonumber
\end{eqnarray}
If we define the function $\bar{f}_\La$ by

\begin{eqnarray}
&&
\bar{f}_\La=
\Delta^{-{N-1\over N^2}}f_\La,
\nonumber
\end{eqnarray}
the KZ equation above is equivalent to 
\begin{eqnarray}
&&
{\partial \bar{f}_\La\over\partial \la_p}
=-{1\over N}\sum_{j\not\in \La_i}{\bar{f}_\La\over\la_p-\la_j}
+{1\over N}\sum_{j\not\in \La_i}{1\over\la_p-\la_j}\bar{f}_{\La^{(pj)}}.
\label{modifkz}
\end{eqnarray}

For a nonnegative integer $r$ and a subset $\{p_1,\cdots,p_r\}\subset\{1,\cdots,Nm\}\ $, set
\begin{eqnarray}
&&
\mdf{\La}{p_1}{\La}{p_r}=\det\Big({\mu_{p_j}^\La(z_i)}\Big)_{1\leq i,j\leq r}.
\nonumber
\end{eqnarray}
Then Theorem \ref{intformula} is equivalent to the following proposition.

\begin{prop}\label{mainprop}
Let
\begin{eqnarray}
&&
\baf_\La=
{\Delta(p_1,\cdots,p_L)^{-1}\over\prod_{i<j}(\La_i\La_j)}
\mdf{\La}{p_1}{\La}{p_L}.
\label{fbar}
\end{eqnarray}
Then, modulo exact forms, $\baf_\La$ satisfies the equation (\ref{modifkz}) and $X\baf_\La=0$ for any $X\in sl_N$.
\end{prop}

Our aim is to prove this proposition.
Let us set
\begin{eqnarray}
&&
\La_r=\{i^r_1,\cdots,i^r_m\}\quad\hbox{ for $1\leq r\leq N$}.
\nonumber
\end{eqnarray}
For the sake of simple exposition we shall prove the equation
(\ref{modifkz}) for $p=i^N_m$. Other cases are similarly proved.

In the following sections we use the usual equality symbol $=$ 
for the equality modulo exact forms.
We remark that all the modulo exact relations, which we use, follow from
the relation in Proposition \ref{exact}.

\subsection{Fundamental Relations}
Now let us give all the relations which we need for our purpose.
For the sake of simplicity, in the formulas below,
we denote $\la_i$ by $i$. 
For instance $i-j=\la_i-\la_j$. 
We set ${\cal K}=\{i^r_l|r<N, (r,l)\neq (N-1,m)\}$ and
\begin{eqnarray}
&&
A_r=
{\prod_{s=1}^m(i^N_m-i^r_s) \over 
\prod_{s=1}^{m-1}(i^N_m-i^N_s)}.
\nonumber
\end{eqnarray}
For $p\neq i^N_m$ and $1\leq j\leq Nm$ we also set
\begin{eqnarray}
&&
\df{p}{j}=\df{\La^{(i^N_mp)}}{j}
\nonumber
\end{eqnarray}
for the sake of simplicity.
Now the relations are given as follows.
\vskip5mm

\noindent
1. For $p\in \La_r$, $r\neq N$
\vskip5mm
\begin{equation}
{\partial\over\partial\la_{i^N_m}}
\df{\La}{p}=
\big(1-{1\over N}\big)
{1\over i^N_m-p}\df{\La}{p}
-{1\over N}{1\over i^N_m-p}
\prod_{j\not\in\La_r,j\neq i^N_m}
{p-j\over i^N_m-j}\df{p}{i^N_m}.
\label{rel1}
\end{equation}

\noindent
2. For $l\neq l'$, $(r,l),(r,l')\neq (N,m)$,
\begin{equation}
\df{i^r_{l'}}{i^r_{l}}
=\!\!
\df{\La}{i^r_l}
+
{i^N_m-i^r_{l'} \over i^N_m-i^{N-1}_{m}}
\sum_{k=1}^{m-1}
{
\prod_{j\not\in\La_r,j\neq i^N_m}(i^r_l-j)
\prod_{s\neq l,l'}^m(i^{N-1}_k-i^r_s) 
\over 
\prod_{j\neq i^N_m,i^{N-1}_m,i^{N-1}_k}(i^{N-1}_k-j)
}
\Big[\!
-\df{\La}{i^{N-1}_k}
\!+\!
\df{i^{N-1}_m}{i^{N-1}_k}
\Big].
\label{rel2}
\end{equation}

\noindent
3. For $l\neq l'$,
$$
\df{i^{N-1}_{l'}}{i^{N-1}_{l}}=
{(i^{N-1}_m-i^{N-1}_{l'})(i^{N}_m-i^{N-1}_{l})
\over
(i^{N-1}_l-i^{N-1}_{l'})(i^{N}_m-i^{N-1}_{m})}
\df{\La}{i^{N-1}_l}
+
{(i^{N-1}_{l}-i^{N-1}_m)(i^{N}_m-i^{N-1}_{l'})
\over
(i^{N-1}_l-i^{N-1}_{l'})(i^{N}_m-i^{N-1}_{m})}
\df{i^{N-1}_m}{i^{N-1}_l}
$$
\begin{equation}
+{(i^{N-1}_{l'}-i^{N-1}_m)(i^N_m-i^{N-1}_{l'})
 \over 
(i^{N-1}_{l'}-i^{N-1}_l)(i^N_m-i^{N-1}_{m})}
\prod_{j\not\in\La_{N-1},j\neq i^N_m}
{i^{N-1}_l-j \over i^{N-1}_{l'}-j}
\Big[
-\df{\La}{i^{N-1}_{l'}}
+
\df{i^{N-1}_{m}}{i^{N-1}_{l'}}
\Big].
\label{rel3}
\end{equation}

\noindent
4. For $r\neq N$,
\begin{eqnarray}
\df{i^{r}_l}{i^N_m}=&&
\sum_{k\in{\cal K}}
\prod_{j\in{\cal K},j\neq k}{i^N_m-j\over k-j}\df{\La}{k}
\nonumber
\\
&&
+
\sum_{k=1}^{m-1}
\prod_{j\in{\cal K},j\neq i^{N-1}_k}
{i^N_m-j \over i^{N-1}_k-j}
B(r,l,k)
\Big[
-\df{\La}{i^{N-1}_k}
+\df{i^{N-1}_m}{i^{N-1}_k}
\Big],
\label{rel4}
\end{eqnarray}
where we set
\begin{eqnarray}
&&
B(r,l,k)=1-
{i^N_m-i^r_l \over A_r}
{\prod_{s\neq l}^{m}(i^{N-1}_k-i^r_s) \over 
\prod_{s=1}^{m-1}(i^{N-1}_k-i^N_s)}.
\nonumber
\end{eqnarray}

\noindent
5. For $(r,l)\neq(N,m)$,
\begin{equation}
\df{i^r_l}{i^r_l}=
\sum_{k\in{\cal K},k\neq i^{r}_l}
{i^N_m-i^r_l \over i^N_m-k}
\prod_{j\in{\cal K},j\neq k,i^{r}_l}
{i^r_l-j \over k-j}
\df{i^r_l}{k}
+
\prod_{j\in{\cal K},j\neq i^{r}_l}
{i^r_l-j \over i^N_m-j}
\df{i^r_l}{i^N_m}.
\label{rel5}
\end{equation}
\vskip5mm

The proofs of these relations are given in appendix.

\subsection{Equation for $\baf_\La$}
Let us take $(p_1,\cdots,p_L)=(i^1_1,\cdots,i^{N-1}_{m-1})$ 
in the expression of $\baf_\La$.
By differentiating the defining equation (\ref{fbar})
in $\la_{i^N_m}$ we have
\begin{eqnarray}
&&
{\partial\baf_\La \over \partial \la_{i^N_m}}=
\sum_{j\not\in\La_N}{-1\over i^N_m-j}\baf_\La
+
{\Delta(i^1_1,\cdots,i^{N-1}_{m-1})^{-1} \over \prod_{t<u}(\La_t\La_u)}
{\partial \over \partial \la_{i^N_m}}
\mdf{\La}{i^1_1}{\La}{i^{N-1}_{m-1}}.
\label{deriv1}
\end{eqnarray}
Substituting (\ref{rel1}) into (\ref{deriv1}) we have
\begin{eqnarray}
{\prod_{t<u}(\La_t\La_u)\over \Delta(i^1_1,\cdots,i^{N-1}_{m-1})^{-1}}
{\partial\baf_\La \over \partial \la_{i^N_m}}
&=&
\Big(
-{1\over i^N_m-i^{N-1}_m}
-{1\over N}\sum_{j\in{\cal K}}{1\over i^N_m-j}
\Big)\mdf{\La}{i^1_1}{\La}{i^{N-1}_{m-1}}
\nonumber
\\
&&
-{1\over N}\sum_{i^r_l\in{\cal K}}{1\over i^N_m-i^r_l}
\prod_{j\not\in\La_r,j\neq i^N_m}
{i^r_l-j \over i^N_m-j}\amdf{i^r_l}{i^N_m}.
\label{deriv2}
\end{eqnarray}
Using (\ref{rel4}) in the second term
of (\ref{deriv2}) we obtain
\vskip3mm
\begin{eqnarray}
&&
(-N)
{\prod_{t<u}(\La_t\La_u)\over \Delta(i^1_1,\cdots,i^{N-1}_{m-1})^{-1}}
{\partial\baf_\La \over \partial \la_{i^N_m}}
\nonumber
\\
&&
=
\Big[
{-m+3 \over i^N_m-i^{N-1}_m}
+\sum_{l=1}^{m-1}{2 \over i^N_m-i^{N-1}_l}
+\sum_{r=1}^{N-2}\sum_{l=1}^m{1\over i^N_m-i^r_l}
\nonumber
\\
&&\quad
+\sum_{r=1}^{N-2}\sum_{l=1}^m
{A_r \over (i^N_m-i^r_l)^2}
{\prod_{s=1}^{m-1} (i^r_l-i^N_s)
 \over \prod_{s\neq l}^{m}(i^r_l-i^r_s)}
\Big]
\mdf{\La}{i^1_1}{\La}{i^{N-1}_{m-1}}
\nonumber
\\
&&
\quad
+
\sum_{i^r_l\in{\cal K}}
{A_r \over (i^N_m-i^r_l)(i^N_m-i^{N-1}_m)}
\sum_{k=1}^{m-1}
{1 \over i^N_m-i^{N-1}_k}
{\prod_{j\not\in\La_r,j\neq i^N_m}(i^r_l-j)
\over 
\prod_{j\in{\cal K},j\neq i^{N-1}_k}
(i^{N-1}_k-j)}B(r,l,k)\times
\nonumber
\\
&&
\quad
\times
\mmdf{\La}{i^{N-1}_m}{i^{N-1}_k}{\La},
\label{calcueq}
\end{eqnarray}
\vskip3mm
where in the second term $\df{i^{N-1}_m}{i^{N-1}_k}$
is on the $l$-th position counted from $\df{\La}{i^r_1}$.
In the derivation of (\ref{calcueq}) we have used
\begin{eqnarray}
&&
\sum_{r=1}^{N-2}\sum_{l=1}^m
{A_r \over (i^N_m-i^r_l)^2}
{i^r_l-i^{N-1}_m \over i^N_m-i^{N-1}_m}
{\prod_{s=1}^{m-1} (i^r_l-i^N_s)
 \over \prod_{s\neq l}^{m}(i^r_l-i^r_s)}
\nonumber
\\
&&
=-{N-2\over i^N_m-i^{N-1}_m}
+
\sum_{r=1}^{N-2}\sum_{l=1}^m
{A_r \over (i^N_m-i^r_l)^2}
{\prod_{s=1}^{m-1} (i^r_l-i^N_s)
 \over \prod_{s\neq l}^{m}(i^r_l-i^r_s)}
\nonumber
\end{eqnarray}
which follows from
\begin{eqnarray}
&&
{i^r_l-i^{N-1}_m \over (i^N_m-i^r_l)(i^N_m-i^{N-1}_m)}
={1\over i^N_m-i^r_l}-{1\over i^N_m-i^{N-1}_m}
\nonumber
\end{eqnarray}
and 
\begin{eqnarray}
&&
\sum_{l=1}^m
{\prod_{s=1}^{m-1} (i^N_m-i^N_s)
 \over (i^N_m-i^r_l)\prod_{s\neq l}^{m}(i^r_l-i^r_s)}
=
{\prod_{s=1}^{m-1} (i^r_l-i^N_s)
 \over \prod_{s=1}^{m}(i^N_m-i^r_s)}
=A_r^{-1}.
\label{resth}
\end{eqnarray}
The equation (\ref{resth}) follows from the residue theorem for the function
$$
{\prod_{s=1}^{m-1}(z-i^N_s)\over (z-i^N_m)\prod_{s=1}^m(z-i^r_s)}.
$$
This is the typical argument to prove an identity
in the proof below.

\subsection{KZ Equation}
We shall rewrite the KZ equation in a similar manner to (\ref{calcueq}).
We take $(p_1,\ldots,p_L)=(i^1_1,\ldots,i^{N-1}_{m-1})$ in the expression
(\ref{fbar}) for $\baf_{\La^{(ji^N_m)}}$ for any $j\notin\La_N$.
Then the KZ equation substituted by (\ref{fbar}) is
\begin{eqnarray}
{\prod_{t<u}(\La_t\La_u)\over \Delta(i^1_1,\cdots,i^{N-1}_{m-1})^{-1}}
{\partial\baf_\La \over \partial \la_{i^N_m}}
&=&
-{1\over N}\sum_{j\not\in\La_N}{1\over i^N_m-j}
\mdf{\La}{i^1_1}{\La}{i^{N-1}_{m-1}}
\nonumber
\\
&&
+{1\over N}\sum_{r=1}^{N-1}\sum_{l=1}^m
{1\over i^N_m-i^r_l}
{\prod_{t<u}(\La_t\La_u) 
\over 
\prod_{t<u}(\La_t^{(i^N_mi^r_l)}\La_u^{(i^N_mi^r_l)})}
\mdf{i^r_l}{i^1_1}{i^r_l}{i^{N-1}_{m-1}}
\label{KZeq1}
\end{eqnarray}
Note the relations
\begin{eqnarray}
&&
{\prod_{t<u}(\La_t\La_u) 
\over 
\prod_{t<u}(\La_t^{(i^N_mi^r_l)}\La_u^{(i^N_mi^r_l)})}
=-
{A_r \over i^N_m-i^r_l}
{\prod_{s=1}^{m-1} (i^r_l-i^N_s)
 \over \prod_{s\neq l}^{m}(i^r_l-i^r_s)}.
\label{ratio}
\end{eqnarray}
Thus we have
\begin{eqnarray}
{\prod_{t<u}(\La_t\La_u)\over \Delta(i^1_1,\cdots,i^{N-1}_{m-1})^{-1}}
{\partial\baf_\La \over \partial \la_{i^N_m}}
&&=
-{1\over N}\sum_{j\not\in\La_N}{1\over i^N_m-j}
\mdf{\La}{i^1_1}{\La}{i^{N-1}_{m-1}}
\nonumber
\\
&&
-{1\over N}\sum_{r=1}^{N-1}\sum_{l=1}^m
{A_r \over (i^N_m-i^r_l)^2}
{\prod_{s=1}^{m-1} (i^r_l-i^N_s)
 \over \prod_{s\neq l}^{m}(i^r_l-i^r_s)}
\mdf{i^r_l}{i^1_1}{i^r_l}{i^{N-1}_{m-1}}.
\label{KZeq2}
\end{eqnarray}

By the equation (\ref{rel5}), 
if $r\neq N$ and $i^r_l\neq i^{N-1}_m$, we have
\begin{eqnarray}
&&
\bmdf{i^r_l}{i^r_l}{i^r_l}{i^r_l}
=\prod_{j\in{\cal K},j\neq i^r_l}
{i^r_l-j \over i^N_m-j}
\bmdf{i^r_l}{i^r_l}{i^N_m}{i^r_l}.
\label{conseq1}
\end{eqnarray}

Using the relation (\ref{rel4}) we have, for $r\leq N-2$,
\begin{eqnarray}
&&
\mmdf{i^r_l}{i^r_l}{i^N_m}{i^r_l}
\nonumber
\\
&&=
\sum_{k=1}^m
\prod_{j\in{\cal K},j\neq i^r_k}
{i^N_m-j \over i^r_k-j}
\mmdf{i^r_l}{\La}{i^r_k}{i^r_l}
\nonumber
\\
&&\quad
+\!
\sum_{k=1}^{m-1}\!
\prod_{j\in{\cal K},j\neq i^{N-1}_k}
{i^N_m-j \over i^{N-1}_k-j}
B(r,l,k)
\mmdf{i^r_l}{i^{N-1}_m}{i^{N-1}_k}{i^r_l},
\label{conseq2}
\end{eqnarray}
and, for $l\neq m$,
\begin{eqnarray}
&&
\smdf{i^{N-1}_l}{i^{N-1}_l}{i^N_m}{i^{N-1}_l}
\nonumber 
\\
=&&\!\!\!\!
\prod_{j\in{\cal K},j\neq i^{N-1}_l}
{i^N_m-j \over i^{N-1}_l-j}
{i^N_m-i^{N-1}_l \over A_{N-1}}
{\prod_{s\neq l}^{m}(i^{N-1}_l-i^{N-1}_s) \over 
\prod_{s=1}^{m-1}(i^{N-1}_l-i^N_s)}
\smdf{i^{N-1}_l}{\La}{i^{N-1}_l}{i^{N-1}_l}
\nonumber
\\
+&&\!\!\!\!\!\!\!\!\!
\sum_{k=1}^{m-1}\!\!
\prod_{j\in{\cal K},j\neq i^{N-1}_k}\!\!
{i^N_m-j \over i^{N-1}_k-j}
B(N-1,l,k)
\smdf{i^{N-1}_l}{i^{N-1}_m}{i^{N-1}_k}{i^{N-1}_l}.
\label{conseq3}
\end{eqnarray}
If we substitute (\ref{conseq1}), (\ref{conseq2}) and (\ref{conseq3})
into (\ref{KZeq2}) we have
\vskip5mm
\begin{eqnarray}
&&
(-N)
{\prod_{t<u}(\La_t\La_u)\over \Delta(i^1_1,\cdots,i^{N-1}_{m-1})^{-1}}
{\partial\baf_\La \over \partial \la_{i^N_m}}
=
\sum_{j\not\in\La_N}{1\over i^N_m-j}
\mdf{\La}{i^1_1}{\La}{i^{N-1}_{m-1}}
\nonumber
\\
&&
+{A_{N-1} \over (i^N_m-i^{N-1}_m)^2}
\prod_{s=1}^{m-1}
{i^{N-1}_m-i^N_s \over
i^{N-1}_m-i^{N-1}_s}
\vsmdf
\nonumber
\\
&&
+\sum_{r=1}^{N-2}\sum_{l=1}^m
{1\over i^N_m-i^r_l}
\sum_{k=1}^m
{A_r \over i^N_m-i^r_k}
{\prod_{j\not\in\La_r,j\neq i^N_m,i^{N-1}_m}
(i^r_l-j) 
\over
\prod_{j\in{\cal K},j\neq i^r_k}(i^r_k-j)}
\mmdf{i^r_l}{\La}{i^r_k}{i^r_l}
\nonumber
\\
&&
+\sum_{l=1}^{m-1}
{1\over i^N_m-i^{N-1}_l}
\smdf{i^{N-1}_l}{\La}{i^{N-1}_l}{i^{N-1}_l}
\nonumber
\\
&&
+\sum_{i^r_l\in{\cal K}}
{1\over (i^N_m-i^r_l)(i^r_l-i^{N-1}_m)}
\sum_{k=1}^{m-1}
{A_r \over i^N_m-i^{N-1}_k}
{\prod_{j\not\in\La_r,j\neq i^N_m}
(i^r_l-j)
\over
\prod_{j\in{\cal K},j\neq i^{N-1}_k}
(i^{N-1}_k-j)}B(r,l,k)
\times
\nonumber
\\
&&
\times
\mmdf{i^r_l}{i^{N-1}_m}{i^{N-1}_k}{i^r_l}
\label{KZeq}
\end{eqnarray}
\vskip5mm
\noindent
where as in the previous case $\df{\La}{i^r_k}$ etc. are
all on the $l$-th place counted from $\df{i^r_l}{i^r_1}$ etc.
For the economy of space we set
\begin{eqnarray}
\Omega(r,l,\La,k)&=&\mmdf{i^r_l}{\La}{i^r_k}{i^r_l},
\nonumber
\\
\Omega(r,l,N-1,k)&=&\mmdf{i^r_l}{i^{N-1}_m}{i^{N-1}_k}{i^r_l},
\nonumber
\end{eqnarray}
where the positions of $\df{\La}{i^r_k}$ and $\df{i^{N-1}_m}{i^{N-1}_k}$
are as above. If $r=N-1$ then $i^r_m$ should be replaced by $i^r_{m-1}$.

The forms
\begin{eqnarray}
&&
\df{\La}{i^1_1},\cdots,\df{\La}{i^{N-1}_{m-1}},
\df{i^{N-1}_m}{i^{N-1}_{1}},\cdots,
\df{i^{N-1}_m}{i^{N-1}_{m-1}}
\nonumber
\end{eqnarray}
are linearly independent in the cohomology group $\mathop{H}^1(C,{\bf C})$.
Since we do not use this fact in this paper, we do not give a proof of it.
But the fact helps to understand the strategy of the proof below.
The right hand side of (\ref{calcueq}) is written using thses forms only.
We shall rewrite the right hand side of (\ref{KZeq}) in terms of these
forms.
\vskip5mm

\subsection{Reduction of Expressions of Fundamental Determinants}
For $1\leq r\leq N-1$, $I=\{s_1<\cdots<s_p\}$, $J=\{t_1<\cdots<t_q\}$,
$p+q=m (r<N-1)$, $p+q=m-1 (r=N-1)$, we set
$$
\Omega^r_{IJ}=
\pdf=
\mmmdf{\La}{i^r_{s_1}}{\La}{i^r_{s_p}}{i^{N-1}_m}{i^{N-1}_{t_1}}{i^{N-1}_m}{i^{N-1}_{t_q}},
$$
where $\df{\La}{i^r_{s_1}}$ is on the $(r-1)m+1$-th position 
counted from $\df{\La}{i^1_1}$.

\vskip3mm
\noindent
$\bullet$
The coefficient of $\Omega^r_{IJ}$ in $\Omega(r,l,\La,k)$.
\vskip3mm
\noindent
We assume $r<N-1$.
Let us denote this coefficient by $\det^1_{IJ}$.
We set
\begin{eqnarray}
&&
F^r_{tu}={1\over i^{N-1}_u-i^r_t},
\quad
G^r_t={1\over \prod_{j\not\in\La_r,j\neq i^N_m}(i^r_t-j)}
\nonumber
\end{eqnarray}
and
$$
A^{rl}_{tu}={i^N_m-i^r_l \over i^N_m-i^{N-1}_m}
{\prod_{s\neq l}^m(i^{N-1}_u-i^r_s) \over
\prod_{j\neq i^N_m,i^{N-1}_m,i^{N-1}_u}(i^{N-1}_u-j)}
{F^r_{tu}\over G^r_t}.
$$
Then (\ref{rel2}) can be written as
\begin{eqnarray}
&&
\df{i^r_l}{i^r_{l'}}=\df{\La}{i^r_{l'}}+
\sum_{k=1}^{m-1}A^{rl}_{tk}
\Big[\!\!\!
-\df{\La}{i^{N-1}_k}
\!\!\!+\!\!\!
\df{i^{N-1}_m}{i^{N-1}_k}
\Big].
\label{rrel2}
\end{eqnarray}

Using this equation we have

\begin{eqnarray}
\hbox{$\det^1_{IJ}$}&=&
\left|
\begin{array}{cccccccccccccc}
1&\quad&\quad&\quad&\quad&\quad&\quad&\quad&0&\vert
&A^{rl}_{11}&\cdots&A^{rl}_{1m-1}\\
\quad&\cdots&\quad&\quad&\quad&\quad&\quad&\quad&\cdot&\vert&
\cdot&\quad&\cdot\\
\quad&\quad&1&\quad&\quad&\quad&\quad&\quad&0&\vert&
\cdot&\quad&\cdot\\
\quad&\quad&\quad&0&\cdot&\cdot&1&\quad&0&\vert&
0&\cdots&0\\
\quad&\quad&\quad&\quad&1&\quad&\vdots&\quad&0&\vert&
\cdot&\quad&\cdot\\
\quad&\quad&\quad&\quad&\quad&\cdots&\vdots&\quad&\cdot&\vert&
\cdot&\quad&\cdot\\
\quad&\quad&\quad&\quad&\quad&\quad&1&\quad&0&\vert&
\cdot&\quad&\cdot\\
\quad&\quad&\quad&\quad&\quad&\quad&\quad&\cdots&\cdot&\vert&
\cdot&\quad&\cdot\\
\quad&\quad&\quad&\quad&\quad&\quad&\quad&\quad&1&\vert&
A^{rl}_{m1}&\cdots&A^{rl}_{mm-1}\\
\end{array}
\right|
\nonumber
\\
&=&
\Big(
{i^N_m-i^r_l \over i^N_m-i^{N-1}_m}
\Big)^{|J|}
\prod_{u\in J}
\Big(
{\prod_{s\neq l}^m(i^{N-1}_u-i^r_s)
\over \prod_{j\neq i^N_m,i^{N-1}_m,i^{N-1}_u}(i^{N-1}_u-j)}
\Big)
{1 \over \prod_{t=1}^mG^r_t}
\nonumber
\\
&&
\times\left|
\begin{array}{cccccccccccccc}
G^r_1&\quad&\quad&\quad&\quad&\quad&\quad&\quad&0&\vert
&F^r_{11}&\cdots&F^r_{1m-1}\\
\quad&\cdots&\quad&\quad&\quad&\quad&\quad&\quad&\cdot&\vert&
\cdot&\quad&\cdot\\
\quad&\quad&G^r_{l-1}&\quad&\quad&\quad&\quad&\quad&0&\vert&
\cdot&\quad&\cdot\\
\quad&\quad&\quad&0&\cdot&\cdot&G^r_l&\quad&0&\vert&
0&\cdots&0\\
\quad&\quad&\quad&\quad&G^r_{l+1}&\quad&\vdots&\quad&0&\vert&
\cdot&\quad&\cdot\\
\quad&\quad&\quad&\quad&\quad&\cdots&\vdots&\quad&\cdot&\vert&
\cdot&\quad&\cdot\\
\quad&\quad&\quad&\quad&\quad&\quad&G^r_k&\quad&0&\vert&
\cdot&\quad&\cdot\\
\quad&\quad&\quad&\quad&\quad&\quad&\quad&\cdots&\cdot&\vert&
\cdot&\quad&\cdot\\
\quad&\quad&\quad&\quad&\quad&\quad&\quad&\quad&G^r_m&\vert&
F^r_{m1}&\cdots&F^r_{mm-1}\\
\end{array}
\right|.
\nonumber
\end{eqnarray}

The meaning of the determinant symbol of the matrix above is
the following.
The matrix consists of two matrices, say, 
the left matrix and the right matrix.
We take $I$-th columns from the left matrix and
$J$-th columns from the right matrix. 
Then form the determinant of
the resulting matrix of degree $|I|+|J|$.
We shall use similar notations from now on.

Notice that, by definition, $\det^1_{IJ}=0$ unless $k\in I$.

\vskip5mm

\noindent
$\bullet$
The coefficient of $\Omega^r_{IJ}$ in $\Omega(r,l,N-1,k)$.
\vskip3mm
\noindent
We assume $r<N-1$.
Let us denote this coefficient by $\det^2_{IJ}$.
Then, again by (\ref{rrel2}),
\begin{eqnarray}
\hbox{$\det^2_{IJ}$}&=&
\Big(
{i^N_m-i^r_l \over i^N_m-i^{N-1}_m}
\Big)^{|J|-1}
\prod_{u\in J\backslash\{k\}}
\Big(
{\prod_{s\neq l}^m(i^{N-1}_u-i^r_s)
\over \prod_{j\neq i^N_m,i^{N-1}_m,i^{N-1}_u}(i^{N-1}_u-j)}
\Big)
{1 \over \prod_{t\neq l}^mG^r_t}
\nonumber
\\
&&
\times\left|
\begin{array}{ccccccccccccc}
G^r_1&\quad&\quad&\cdot&\quad&\quad&\quad
&\vert&F^{rl}_{11}&\cdots&\bar{A}^{rl}_{1k}&\cdots&F^{rl}_{1m-1}\\
\quad&\cdots&\quad&\cdot&\quad&\quad&\quad
&\vert&\cdot&\quad&\cdot&\quad&\cdot\\
\quad&\quad&G^r_{l-1}&\cdot&\quad&\quad&\quad
&\vert&\cdot&\quad&\cdot&\quad&\cdot\\
\cdot&\cdot&\cdot&0&\cdot&\cdot&\cdot
&\vert&0&\quad&1&\quad&0\\
\quad&\quad&\quad&\quad&G^r_{l+1}&\quad&\quad
&\vert&\cdot&\quad&\cdot&\quad&\cdot\\
\quad&\quad&\quad&\quad&\quad&\cdots&\quad
&\vert&\cdot&\quad&\cdot&\quad&\cdot\\
\quad&\quad&\quad&\quad&\quad&\quad&G^r_m
&\vert&F^{rl}_{m1}&\cdots&\bar{A}^{rl}_{mk}&\cdots&F^{rl}_{mm-1}\\
\end{array}
\right|,
\nonumber
\end{eqnarray}
where
\begin{eqnarray}
&&
\bar{A}^{rl}_{tk}
=A^{rl}_{tk}G^r_t.
\nonumber
\end{eqnarray}
\vskip5mm

\noindent
$\bullet$
The coefficient of $\Omega^{N-1}_{IJ}$ in 
$\Omega(N-1,l,\La,l)$.
\vskip3mm
\noindent
This coefficient is denoted by $\det^3_{IJ}$.
Let us set, for $s\neq l$,

\begin{eqnarray}
&&
C^l_{ss}=
{(i^{N-1}_l-i^{N-1}_m)(i^N_m-i^{N-1}_s)
\over
(i^N_m-i^{N-1}_m)(i^{N-1}_l-i^{N-1}_s)},
\quad
C^l_{sl}=
-{(i^N_m-i^{N-1}_l)(i^{N-1}_l-i^{N-1}_m)
\over
(i^N_m-i^{N-1}_m)(i^{N-1}_l-i^{N-1}_s)}
\prod_{j\not\in\La_{N-1},j\neq i^N_m}
{i^{N-1}_s-j \over i^{N-1}_l-j},
\nonumber
\\
&&
D^l_{ss}=
{(i^N_m-i^{N-1}_l)(i^{N-1}_s-i^{N-1}_m)
\over
(i^N_m-i^{N-1}_m)(i^{N-1}_s-i^{N-1}_l)},
\quad
D^l_{sl}=-C^l_{sl}.
\nonumber
\end{eqnarray}

Then, (\ref{rel3}) is written as
\begin{eqnarray}
&&
\df{i^{N-1}_l}{i^{N-1}_s}=
C^l_{ss}\df{\La}{i^{N-1}_s}+C^l_{sl}\df{\La}{i^{N-1}_l}
+D^l_{ss}\df{i^{N-1}_m}{i^{N-1}_s}+D^l_{sl}\df{i^{N-1}_m}{i^{N-1}_l}.
\label{rrel3}
\end{eqnarray}

Using (\ref{rrel3}) we have
\begin{eqnarray}
&&
\hbox{$\det^3_{IJ}$}=
\nonumber
\\
&&
\left|
\begin{array}{ccccccccccccc}
C^l_{11}&\quad&\quad&C^l_{1l}&\quad&\quad&0
&\vert&D^l_{11}&\cdots&D^l_{1l}&\cdots&0\\
\quad&\cdots&\quad&\cdot&\quad&\quad&\cdot
&\vert&\cdot&\quad&\cdot&\quad&\cdot\\
\quad&\quad&C^l_{l-1l-1}&\cdot&\quad&\quad&\cdot
&\vert&\cdot&\quad&\cdot&\quad&\cdot\\
0&\cdot&\cdot&1&\cdot&\cdot&0
&\vert&0&\quad&0&\quad&0\\
\cdot&\quad&\quad&\cdot&C^l_{l+1l+1}&\quad&\cdot
&\vert&\cdot&\quad&\cdot&\quad&\cdot\\
\cdot&\quad&\quad&\cdot&\quad&\cdots&\cdot
&\vert&\cdot&\quad&\cdot&\quad&\cdot\\
0&\quad&\quad&C^l_{m-1l}&\quad&\quad&C^l_{m-1m-1}
&\vert&0&\cdots&D^l_{m-1l}&\cdots&D^l_{m-1m-1}\\
\end{array}
\right|.
\nonumber
\end{eqnarray}

\noindent
$\bullet$
The coefficient of $\Omega^{N-1}_{IJ}$ in 
$\Omega(N-1,l,N-1,k)$.
\vskip3mm
\noindent
This coefficient is denoted by $\det^4_{IJ}$.
Then
\begin{eqnarray}
&&
\hbox{$\det^4_{IJ}$}=
\nonumber
\\
&&
\left|
\begin{array}{cccccccccccccccc}
C^l_{11}&\quad&\quad&C^l_{1l}&\quad&\quad&0
&\vert&D^l_{11}&\cdot&0&\cdot&D^l_{1l}&\cdots&0\\
\quad&\cdots&\quad&\cdot&\quad&\quad&\cdot
&\vert&\cdot&\cdots&\cdot&\quad&\cdot&\quad&\cdot\\
\quad&\quad&C^l_{l-1l-1}&\cdot&\quad&\quad&\cdot
&\vert&\cdot&\quad&\cdot&\cdots&\cdot&\quad&\cdot\\
0&\cdot&\cdot&0&\cdot&\cdot&0
&\vert&0&\cdots&1&\cdots&0&\quad&0\\
\cdot&\quad&\quad&\cdot&C^l_{l+1l+1}&\quad&\cdot
&\vert&\cdot&\quad&\cdot&\quad&\cdot&\cdots&\cdot\\
\cdot&\quad&\quad&\cdot&\quad&\cdots&\cdot
&\vert&\cdot&\quad&\cdot&\quad&\cdot&\cdots&\cdot\\
0&\quad&\quad&C^l_{m-1l}&\quad&\quad&C^l_{m-1m-1}
&\vert&0&\cdots&0&\cdots&D^l_{m-1l}&\cdots&D^l_{m-1m-1}\\
\end{array}
\right|.
\nonumber
\end{eqnarray}
In the right matrix $1$ is in the $(l,k)$ component.

\subsection{Comparison of Two Equations}
Now let us calculate the reduced expression of the right hand
side of (\ref{KZeq}) and compare it with (\ref{calcueq}).
We shall calculate the coefficient of $\Omega^r_{IJ}$ by dividing the
case into nine as
\begin{description}
\item[(I)] $r=N-1$, $I=\phi$, $J=\{1,2,\cdots,m-1\}$,
\item[(II)] $r=N-1$, $I=\{1,2,\cdots,m-1\}$, $J=\phi$ or
$r<N-1$, $I=\{1,2,\cdots,m\}$, $J=\phi$,
\item[(III)] $r=N-1$, $I=\{1,2,\cdots,m-1\}\backslash\{t\}$, $J=\{t\}$,
\item[(IV)] $r=N-1$, $I=\{1,2,\cdots,m-1\}\backslash\{u\}$, 
$J=\{t\}$, $u\neq t$,
\item[(V)] $r<N-1$, $I=\{1,2,\cdots,m-1\}\backslash\{u\}$, $J=\{t\}$, 
\item[(VI)] $r=N-1$, $I\cap J=\phi$, $|I|\geq 1$, $|J|\geq 2$,
\item[(VII)] $r=N-1$, $|I\cap J|=1$, $|J|\geq 2$,
\item[(VIII)] $r=N-1$, $|I\cap J|\geq 2$,
\item[(IX)] $r<N-1$, $|J|\geq 2$.
\end{description}
\vskip7mm

\noindent
{\bf(I)}: The coefficient of $\Omega^{N-1}_{IJ}$ with 
$I=\phi$, $J=\{1,2,\cdots,m-1\}$.
\vskip3mm
\noindent
Let us calculate the contribution from the term which contains
$\Omega(N-1,l,N-1,k)$ $(1\leq l\leq m-1)$.
We have, for $k\neq l$,
\begin{eqnarray}
\hbox{$\det^4_{IJ}$}&=&
-{D^l_{kl} \over D^l_{kk}}\prod_{s\neq l}^{m-1}D^l_{ss}
\nonumber
\\
&=&
{i^{N-1}_l-i^{N-1}_m \over i^{N-1}_k-i^{N-1}_m}
\Big({i^N_m-i^{N-1}_l \over i^N_m-i^{N-1}_m}\Big)^{m-2}
\prod_{s\neq l}^{m-1}{i^{N-1}_m-i^{N-1}_s \over i^{N-1}_l-i^{N-1}_s}
\prod_{j\not\in\La_{N-1},j\neq i^N_m}{i^{N-1}_k-j\over i^{N-1}_l-j}.
\nonumber
\end{eqnarray}
For $k=l$, $\det^4_{IJ}=\prod_{s\neq l}^{m-1}D^l_{ss}$.
If we set $k=l$ in the right hand side of
\begin{eqnarray}
&&
{D^l_{kl}\over D^l_{kk}}=-
{i^{N-1}_l-i^{N-1}_m \over i^{N-1}_k-i^{N-1}_m}
\prod_{j\not\in\La_{N-1},j\neq i^N_m}{i^{N-1}_k-j\over i^{N-1}_l-j},
\nonumber
\end{eqnarray}
we have $-1$. Thus the formula for $\det^4_{IJ}$ given above
is valid for all $1\leq k,l\leq m-1$.
The contribution to (\ref{KZeq}) is
\begin{eqnarray}
&&
-A_{N-1}
{\prod_{s=1}^{m-1}(i^{N-1}_m-i^{N-1}_s)
\over (i^N_m-i^{N-1}_m)^{m-2}}
\sum_{l=1}^{m-1}
{(i^N_m-i^{N-1}_l)^{m-3} 
\over \prod_{s\neq l}^m(i^{N-1}_l-i^{N-1}_s)}
\times
\nonumber
\\
&&
\times
\sum_{k=1}^{m-1}
{\prod_{s=1}^{m-1}(i^{N-1}_k-i^N_s) 
\over (i^N_m-i^{N-1}_k)\prod_{s\neq k}^m(i^{N-1}_k-i^{N-1}_s)}
B(N-1,l,k).
\label{coef1}
\end{eqnarray}
By the residue theorem we have
\begin{eqnarray}
&&
\sum_{k=1}^{m-1}
{\prod_{s=1}^{m-1}(i^{N-1}_k-i^N_s) 
\over (i^N_m-i^{N-1}_k)\prod_{s\neq k}^m(i^{N-1}_k-i^{N-1}_s)}
={1\over A_{N-1}}
-
{\prod_{s=1}^{m-1}(i^{N-1}_m-i^N_s)
\over
(i^N_m-i^{N-1}_m)\prod_{s=1}^{m-1}(i^{N-1}_m-i^{N-1}_s)},
\nonumber
\\
&&
{i^N_m-i^{N-1}_l \over A_{N-1}}
\sum_{k=1}^{m-1}
{\prod_{s\neq l}(i^{N-1}_k-i^{N-1}_s) 
\over (i^N_m-i^{N-1}_k)\prod_{s\neq k}^m(i^{N-1}_k-i^{N-1}_s)}
={1\over A_{N-1}}.
\nonumber
\end{eqnarray}

Thus, substituting the definition of $B(N-1,l,k)$ to (\ref{coef1}),
\begin{eqnarray}
(\ref{coef1})&=&
A_{N-1}
{\prod_{s=1}^{m-1}(i^{N-1}_m-i^N_s)
\over (i^N_m-i^{N-1}_m)^{m-1}}
\sum_{l=1}^{m-1}
{(i^N_m-i^{N-1}_l)^{m-3} 
\over \prod_{s\neq l}^m(i^{N-1}_l-i^{N-1}_s)}
\nonumber
\\
&=&
-{A_{N-1} \over (i^N_m-i^{N-1}_m)^2}
\prod_{s=1}^{m-1}{i^{N-1}_m-i^N_s \over i^{N-1}_m-i^{N-1}_s}.
\label{coef2}
\end{eqnarray}
Here we again use the residue theorem to evaluate the summation
in $l$. Hence the coefficient of $\Omega$ in the right hand side
of (\ref{KZeq}) is zero. This is the case for (\ref{calcueq}).
\vskip7mm

\noindent
{\bf (II)}: The coefficient of 
$\Omega^r_{IJ}=\mdf{\La}{i^1_1}{\La}{i^{N-1}_{m-1}}$ with $r=N-1$
$I=\{1,2,...,m-1\}$ and $J=\phi$ or $r<N-1$
$I=\{1,2,...,m\}$ and $J=\phi$.
\vskip3mm
\noindent
(II-I)
 The contribution to (\ref{KZeq}) from the term which contains
$\Omega(r,l,\La,k)$ $(1\leq l\leq m)$.\par
We have
$$
\hbox{$\det^1_{IJ}$}=
\left\{
\begin{array}{ll}
0 & \mbox{$k\neq l$}\\
1 & \mbox{$k=l$}
\end{array}
\right.
$$
Hence the contribution to the rhs of (\ref{KZeq}) from the term 
containing these determinants is
\begin{eqnarray}
&&
\sum_{r=1}^{N-2}\sum_{l=1}^m
{A_r \over (i^N_m-i^r_l)^2}
{\prod_{j\not\in\La_r,j\neq i^N_m,i^{N-1}_m}(i^r_l-j)
\over \prod_{j\in{\cal K},j\neq i^r_l}(i^r_l-j)}
\nonumber
\\
&&
=\sum_{r=1}^{N-2}\sum_{l=1}^m
{A_r \over (i^N_m-i^r_l)^2}
{\prod_{s=1}^{m-1}(i^r_l-i^N_s)
\over \prod_{s\neq l}^m(i^r_l-i^r_s)}.
\label{cont21}
\end{eqnarray}
\par

\noindent
(II-II)
The contribution to (\ref{KZeq}) from the term which contains
$\Omega(N-1,l,\La,l)$ $(1\leq l\leq m-1)$.\par
We have
\begin{eqnarray}
&&
\hbox{$\det^3_{IJ}$}=
\prod_{s\neq l}^{m-1}C^l_{ss}
=\Big(
{i^{N-1}_l-i^{N-1}_m \over i^N_m-i^{N-1}_m}
\Big)^{m-2}
\prod_{s\neq l}^{m-1}
{i^N_m-i^{N-1}_s \over i^{N-1}_l-i^{N-1}_s}.
\nonumber
\end{eqnarray}
The contribution to the rhs of (\ref{KZeq}) from 
the terms containing these determinants is
\begin{eqnarray}
&&
{\prod_{s=1}^{m-1}
(i^N_m-i^{N-1}_s) \over (i^N_m-i^{N-1}_m)^{m-2}}
\sum_{l=1}^{m-1}
{(i^{N-1}_l-i^{N-1}_m)^{m-2} \over
(i^N_m-i^{N-1}_l)^2
\prod_{s\neq l}^{m-1}(i^{N-1}_l-i^{N-1}_s)}
\nonumber
\\
&&=
-{\prod_{s=1}^{m-1}
(i^N_m-i^{N-1}_s) \over (i^N_m-i^{N-1}_m)^{m-2}}
\mathop{Res}_{z=i^N_m}
{(z-i^{N-1}_m)^{m-2} \over
(z-i^N_m)^2
\prod_{s=1}^{m-1}(z-i^{N-1}_s)}
\nonumber
\\
&&=
-{m-2 \over i^N_m-i^{N-1}_m}
+\sum_{s=1}^{m-1}{1 \over i^N_m-i^{N-1}_s}.
\label{cont22}
\end{eqnarray}
\vskip2mm

From (\ref{cont21}) and (\ref{cont22}) the coefficient
of $\Omega$ in the rhs of (\ref{KZeq}) is
\begin{eqnarray}
&&
\sum_{j\not\in\La_N}{1\over i^N_m-j}
-{m-2 \over i^N_m-i^{N-1}_m}
+\sum_{s=1}^{m-1}{1 \over i^N_m-i^{N-1}_s}
+\sum_{r=1}^{N-2}\sum_{l=1}^m
{A_r \over (i^N_m-i^r_l)^2}
{\prod_{s=1}^{m-1}(i^r_l-i^N_s)
\over \prod_{s\neq l}^m(i^r_l-i^r_s)}
\nonumber
\\
&&=
{-m+3 \over i^N_m-i^{N-1}_m}
+\sum_{s=1}^{m-1}{2 \over i^N_m-i^{N-1}_s}
+\sum_{r=1}^{N-2}\sum_{l=1}^m{1\over i^N_m-i^r_l}
+\sum_{r=1}^{N-2}\sum_{l=1}^m
{A_r \over (i^N_m-i^r_l)^2}
{\prod_{s=1}^{m-1}(i^r_l-i^N_s)
\over \prod_{s\neq l}^m(i^r_l-i^r_s)}.
\nonumber
\end{eqnarray}
This coinsides with the coefficient of $\Omega$ in the
rhs of (\ref{calcueq})
\vskip7mm

\noindent
{\bf(III)}: The coefficient of 
$\Omega^{N-1}_{IJ}$ for which $I=\{1,2,...,m-1\}\backslash\{t\}$
and $J=\{t\}$.
\vskip3mm
\noindent
(III-I)
The contribution to (\ref{KZeq}) 
from the term which contains
$\Omega(N-1,l,\La,l)$.\par
It is obvious that $\det^3_{IJ}=0$ for $t=l$. 
We have, for $t\neq l$,
\begin{eqnarray}
\hbox{$\det^3_{IJ}$}&=&
(-1)^{m+t+1}D^l_{tt}\prod_{s\neq t,l}^{m-1}C^l_{ss}
\nonumber
\\
&=&
(-1)^{m+t+1}
{(i^{N-1}_m-i^{N-1}_t)(i^{N-1}_l-i^{N-1}_m)^{m-3}
\over
(i^N_m-i^{N-1}_m)^{m-2}}
{\prod_{s\neq t}^{m-1}(i^N_m-i^{N-1}_s)
\over
\prod_{s\neq l}^{m-1}(i^{N-1}_l-i^{N-1}_s)}.
\nonumber
\end{eqnarray}
The contribution to the rhs of (\ref{KZeq}) from 
the terms containing these determinants is
\begin{eqnarray}
&&
(-1)^{m+t+1}
{(i^{N-1}_m-i^{N-1}_t)\prod_{s\neq t}^{m-1}(i^N_m-i^{N-1}_s)
\over
(i^N_m-i^{N-1}_m)^{m-2}}
\sum_{l\neq t}^{m-1}
{(i^{N-1}_l-i^{N-1}_m)^{m-3}
\over
(i^N_m-i^{N-1}_l)\prod_{s\neq l}^{m-1}(i^{N-1}_l-i^{N-1}_s)}.
\label{cont31}
\end{eqnarray}
\vskip2mm

\noindent
(III-II)
The contribution from the term which contains
$\Omega(N-1,l,N-1,l)$ 
$(1\leq l\leq m-1)$.\par
If $l\neq t$ then $\det^4_{IJ}=0$. 
We have, for $l=t$,
\begin{eqnarray}
&&
\hbox{$\det^4_{IJ}$}=
(-1)^{m+t+1}\prod_{s\neq t}^{m-1}C^t_{ss}=
(-1)^{m+t+1}
\Big({i^{N-1}_t-i^{N-1}_m \over i^N_m-i^{N-1}_m}\Big)^{m-2}
\prod_{s\neq t}^{m-1}{i^N_m-i^{N-1}_s \over i^{N-1}_t-i^{N-1}_s}.
\nonumber
\end{eqnarray}
The contribution to the rhs of (\ref{KZeq}) from 
the terms containing these determinants is
\begin{eqnarray}
&&
{(-1)^{m+t+1}A_{N-1}(i^{N-1}_t-i^{N-1}_m)^{m-3}
\over (i^{N-1}_t-i^N_m)^2(i^N_m-i^{N-1}_m)^{m-2}}
{\prod_{s=1}^{m-1}(i^{N-1}_t-i^N_s)
\over
\prod_{s\neq t}^{m-1}(i^{N-1}_t-i^{N-1}_s)}
\prod_{s\neq t}^{m-1}
{i^N_m-i^{N-1}_s \over i^{N-1}_t-i^{N-1}_s}
\nonumber
\\
&&
+
{(-1)^{m+t+1} \over i^{N-1}_t-i^N_m}
\Big({i^{N-1}_t-i^{N-1}_m \over i^N_m-i^{N-1}_m}\Big)^{m-2}
\prod_{s\neq t}^{m-1}
{i^N_m-i^{N-1}_s \over i^{N-1}_t-i^{N-1}_s}.
\label{cont32}
\end{eqnarray}
Note that the second term of (\ref{cont32}) is the $l=t$ case
of the summand of (\ref{cont31}).
\vskip2mm

\noindent
(III-III)
The contribution from the term which contains
$\Omega(N-1,l,N-1,k)$ 
with $1\leq k,l\leq m-1$ and $k\neq l$.\par
If $k\neq t$ then $\det^4_{IJ}=0$. 
We have, for $k=t$,
\begin{eqnarray}
\hbox{$\det^4_{IJ}$}&=&
(-1)^{m+t}C^l_{tl}\prod_{s\neq t,l}^{m-1}C^l_{ss}
\nonumber
\\
&=&
(-1)^{m+t+1}
\Big({i^{N-1}_l-i^{N-1}_m \over i^N_m-i^{N-1}_m}\Big)^{m-2}
{\prod_{s\neq t}^{m-1}(i^N_m-i^{N-1}_s)
\over 
\prod_{s\neq l}^{m-1}(i^{N-1}_l-i^{N-1}_s)}
\prod_{j\not\in\La_{N-1},j\neq i^N_m}
{i^{N-1}_t-j \over i^{N-1}_l-j}.
\nonumber
\end{eqnarray}
The contribution to the rhs of (\ref{KZeq}) is
\begin{eqnarray}
&&
{(-1)^{m+t+1}A_{N-1}
\prod_{s\neq t}^{m-1}(i^N_m-i^{N-1}_s)
\prod_{s=1}^{m-1}(i^{N-1}_t-i^N_s)
\over 
(i^N_m-i^{N-1}_t)(i^N_m-i^{N-1}_m)^{m-2}
\prod_{s\neq t}^{m-1}(i^{N-1}_t-i^{N-1}_s)}
\times
\nonumber
\\
&&
\times
\sum_{l\neq t}^{m-1}
{(i^{N-1}_l-i^{N-1}_m)^{m-3}
\over 
(i^N_m-i^{N-1}_l)\prod_{s\neq l}^{m-1}(i^{N-1}_l-i^{N-1}_s)}.
\label{cont33}
\end{eqnarray}
In deriving (\ref{cont33}) we use
\begin{eqnarray}
&&
{i^N_m-i^{N-1}_l \over A_{N-1}}
{\prod_{s\neq l}^{m-1}(i^{N-1}_t-i^{N-1}_s)
\over 
\prod_{s=1}^{m-1}(i^{N-1}_t-i^N_s)}=0
\nonumber
\end{eqnarray}
which is a consequence of $t\neq l$.
Note that the first term in (\ref{cont32}) is the $l=t$ case
of the summand of (\ref{cont33}).
\vskip2mm

We add (\ref{cont31}), (\ref{cont32}),(\ref{cont33}) and obtain
\begin{eqnarray}
&&
{(-1)^{m+t+1}A_{N-1} 
\prod_{s=1}^{m-1}(i^{N-1}_t-i^N_s)
\over (i^N_m-i^{N-1}_t)(i^N_m-i^{N-1}_m)^{m-2}
}
\prod_{s\neq t}^{m-1}
{(i^N_m-i^{N-1}_s)\over(i^{N-1}_t-i^{N-1}_s)}
B(N-1,t,t)
\nonumber
\\
&&
\times
\sum_{l=1}^{m-1}
{(i^{N-1}_l-i^{N-1}_m)^{m-3}
\over
(i^N_m-i^{N-1}_l)\prod_{s\neq l}^{m-1}(i^{N-1}_l-i^{N-1}_s)}
\nonumber
\\
&&=
{(-1)^{m+t+1}
A_{N-1} \over (i^N_m-i^{N-1}_t)^2(i^N_m-i^{N-1}_m)}
{\prod_{s=1}^{m-1}(i^{N-1}_t-i^N_s)
\over
\prod_{s\neq t}^{m-1}(i^{N-1}_t-i^{N-1}_s)}B(N-1,t,t).
\label{cont3}
\end{eqnarray}
In deriving (\ref{cont3}) we use the identity
\begin{eqnarray}
&&
\sum_{l=1}^{m-1}
{(i^{N-1}_l-i^{N-1}_m)^{m-3}
\over 
(i^N_m-i^{N-1}_l)\prod_{s\neq l}^{m-1}(i^{N-1}_l-i^{N-1}_s)}
={(i^N_m-i^{N-1}_m)^{m-3}
\over
\prod_{s=1}^{m-1}(i^N_m-i^{N-1}_s)}.
\label{resth1}
\end{eqnarray}
The equation (\ref{cont3}) is nothing but the corresponding
coefficient in the rhs of (\ref{calcueq}).
\vskip7mm

\noindent
{\bf(IV)}: The coefficient of 
$\Omega^{N-1}_{IJ}$ with $I=\{1,2,...,m-1\}\backslash\{u\}$
and $J=\{t\}$, $u\neq t$.
\vskip3mm
\noindent
(IV-I)
The contribution to (\ref{KZeq}) 
from the term which contains
$\Omega(N-1,l,\La,l)$.\par
If $l=u$ then $\det^3_{IJ}=0$. 
We assume $l\neq u$.
Then $\det^3_{IJ}=0$ for $t\neq l$, since $u$-th and $l$-th rows are
proportional. Thus we assume $l=t$.
We have
\begin{eqnarray}
&&
\hbox{$\det^3_{IJ}$}=(-1)^{u+m+1}D^t_{ut}\prod^{m-1}_{s\neq u,t}C^t_{ss}
\nonumber
\\
&&=
(-1)^{u+m+1}
\Big({i^{N-1}_t-i^{N-1}_m \over i^N_m-i^{N-1}_m}\Big)^{m-2}
{\prod_{s\neq u}^{m-1}(i^N_m-i^{N-1}_s)
\over
\prod_{s\neq t}^{m-1}(i^{N-1}_t-i^{N-1}_s)}
\prod_{j\not\in\La_{N-1},j\neq i^N_m}
{i^{N-1}_u-j \over i^{N-1}_t-j}.
\nonumber
\end{eqnarray}
The contribution to the rhs of (\ref{KZeq}) from 
the terms containing these determinants is
\begin{eqnarray}
&&
{(-1)^{u+m+1} \over i^N_m-i^{N-1}_t}
\Big({i^{N-1}_t-i^{N-1}_m \over i^N_m-i^{N-1}_m}\Big)^{m-2}
{\prod_{s\neq u}^{m-1}(i^N_m-i^{N-1}_s)
\over
\prod_{s\neq t}^{m-1}(i^{N-1}_t-i^{N-1}_s)}
\prod_{j\not\in\La_{N-1},j\neq i^N_m}
{i^{N-1}_u-j \over i^{N-1}_t-j}.
\label{cont41}
\end{eqnarray}
\vskip2mm

\noindent
(IV-II)
The contribution to (\ref{KZeq}) 
from the term which contains
$\Omega(N-1,l,N-1,l)$ with
$1\leq l\leq m-1$.\par
If $l\neq t$ then $\det^4_{IJ}=0$. 
We have, for $l=t$,
\begin{eqnarray}
\hbox{$\det^4_{IJ}$}&=&(-1)^{u+m}C^t_{ut}\prod_{s\neq u,t}^{m-1}C^t_{ss}
=(-1)^{u+m+1}D^t_{ut}\prod_{s\neq u,t}^{m-1}C^t_{ss}
\nonumber
\\
&=&
(-1)^{u+m+1}
\Big({i^{N-1}_t-i^{N-1}_m \over i^N_m-i^{N-1}_m}\Big)^{m-2}
{\prod_{s\neq u}^{m-1}(i^N_m-i^{N-1}_s)
\over
\prod_{s\neq t}^{m-1}(i^{N-1}_t-i^{N-1}_s)}
\prod_{j\not\in\La_{N-1},j\neq i^N_m}
{i^{N-1}_u-j \over i^{N-1}_t-j}.
\nonumber
\end{eqnarray}
The contribution to the rhs of (\ref{KZeq}) from 
the terms containing these determinants is
\begin{eqnarray}
&&
{(-1)^{u+m+1}A_{N-1}(i^{N-1}_t-i^{N-1}_m)^{m-3}
\over
(i^N_m-i^{N-1}_t)^2(i^N_m-i^{N-1}_m)^{m-2}}
{\prod_{s\neq u}^{m-1}(i^N_m-i^{N-1}_s)
\prod_{j\not\in\La_{N-1},j\neq i^N_m}(i^{N-1}_u-j)
\over
\prod_{s\neq t}^{m-1}(i^{N-1}_t-i^{N-1}_s)
\prod_{j\not\in\La_N,j\neq i^{N-1}_t,i^{N-1}_m}(i^{N-1}_t-j)}
\nonumber
\\
&&
+
{(-1)^{u+m} \over i^N_m-i^{N-1}_t}
\Big({i^{N-1}_t-i^{N-1}_m \over i^N_m-i^{N-1}_m}\Big)^{m-2}
{\prod_{s\neq u}^{m-1}(i^N_m-i^{N-1}_s)
\over
\prod_{s\neq t}^{m-1}(i^{N-1}_t-i^{N-1}_s)}
\prod_{j\not\in\La_{N-1},j\neq i^N_m}
{i^{N-1}_u-j \over i^{N-1}_t-j}.
\label{cont42}
\end{eqnarray}
Note that the second term of this equation is the minus
of (\ref{cont41}).
\vskip2mm

\noindent
(IV-III)
The contribution to (\ref{KZeq}) 
from the term which contains
$\Omega(N-1,l,N-1,k)$ with
$1\leq k,l\leq m-1$ and $k\neq l$.\par
If $k\neq t$ then $\det^4_{IJ}=0$. 
We have, for $k=t$,
\begin{eqnarray}
\hbox{$\det^4_{IJ}$}&=&
(-1)^{u+m}{C^l_{ul}\over C^l_{uu}}\prod_{s\neq l}^{m-1}C^l_{ss}
\nonumber
\\
&=&
(-1)^{u+m+1}
\Big({i^{N-1}_l-i^{N-1}_m \over i^N_m-i^{N-1}_m}\Big)^{m-2}
{\prod_{s\neq u}^{m-1}(i^N_m-i^{N-1}_s)
\over
\prod_{s\neq l}^{m-1}(i^{N-1}_l-i^{N-1}_s)}
\prod_{j\not\in\La_{N-1},j\neq i^N_m}
{i^{N-1}_u-j \over i^{N-1}_l-j}.
\nonumber
\end{eqnarray}
Here we understand $C^l_{ul}/C^l_{uu}=-1$ for $u=l$.
This follows from the equation
\begin{eqnarray}
&&
{C^l_{ul} \over C^l_{uu}}=
-{i^N_m-i^{N-1}_l \over i^N_m-i^{N-1}_u}
\prod_{j\not\in\La_{N-1},j\neq i^N_m}
{i^{N-1}_u-j \over i^{N-1}_l-j}.
\nonumber
\end{eqnarray}
The contribution to the rhs of (\ref{KZeq}) is
\begin{eqnarray}
&&
{(-1)^{u+m+1}A_{N-1}\prod_{s\neq u}^{m-1}(i^N_m-i^{N-1}_s)
\prod_{j\not\in\La_{N-1},j\neq i^N_m}(i^{N-1}_u-j)
\over
(i^N_m-i^{N-1}_t)(i^N_m-i^{N-1}_m)^{m-2}
\prod_{j\in{\cal K},j\neq i^{N-1}_t}(i^{N-1}_t-j)}
\nonumber
\\
&&
\times
\sum_{l\neq t}^{m-1}
{(i^{N-1}_l-i^{N-1}_m)^{m-3}
\over 
(i^N_m-i^{N-1}_l)\prod_{s\neq l}^{m-1}(i^{N-1}_l-i^{N-1}_s)}.
\label{cont43}
\end{eqnarray}
Note that the $l=t$ term of this equation is precisely
the first term of (\ref{cont42}).
\vskip2mm

Thus, using (\ref{resth1}), we have
\begin{eqnarray}
&&
\hbox{(\ref{cont41})}+\hbox{(\ref{cont42})}+\hbox{(\ref{cont43})}
\nonumber
\\
&&
=
{(-1)^{u+m+1}A_{N-1}\prod_{s\neq u}^{m-1}(i^N_m-i^{N-1}_s)
\prod_{j\not\in\La_{N-1},j\neq i^N_m}(i^{N-1}_u-j)
\over
(i^N_m-i^{N-1}_t)(i^N_m-i^{N-1}_m)^{m-2}
\prod_{j\in{\cal K},j\neq i^{N-1}_t}(i^{N-1}_t-j)}
\nonumber
\\
&&
\times
\sum_{l=1}^{m-1}
{(i^{N-1}_l-i^{N-1}_m)^{m-3}
\over 
(i^N_m-i^{N-1}_l)\prod_{s\neq l}^{m-1}(i^{N-1}_l-i^{N-1}_s)}
\nonumber
\\
&&
=
{(-1)^{u+m+1}A_{N-1} 
\over
(i^N_m-i^{N-1}_u)(i^N_m-i^{N-1}_m)(i^N_m-i^{N-1}_t)}
{\prod_{j\not\in\La_{N-1},j\neq i^N_m}(i^{N-1}_u-j)
\over
\prod_{j\in{\cal K},j\neq i^{N-1}_t}(i^{N-1}_t-j)}
\nonumber
\end{eqnarray}
which coinsides with the coefficient of $\Omega^{N-1}_{IJ}$
in the rhs of (\ref{calcueq}).
\vskip7mm

\noindent
{\bf(V)}: The coefficient of 
$\Omega^r_{IJ}$ for which $I=\{1,2,...,m\}\backslash\{u\}$
and $J=\{t\}$, $r<N-1$.
\vskip3mm
\noindent
(V-I)
The contribution from the term which contains
$\Omega(r,l,\La,l)$.\par
If $l=u$ then $\det^1_{IJ}=0$. 
We assume $l\neq u$.
Then
\begin{eqnarray}
\hbox{$\det^1_{IJ}$}&=&
(-1)^{u+m}
{i^N_m-i^r_l \over i^N_m-i^{N-1}_m}
{\prod_{s\neq l}^m(i^{N-1}_t-i^r_s)
\over
\prod_{j\neq i^N_m,i^{N-1}_m,i^{N-1}_t}(i^{N-1}_t-j)}
{1\over \prod_{s=1}^mG^r_s}
F^r_{ut}\prod_{s\neq u}^mG^r_s
\nonumber
\\
&=&
(-1)^{u+m}
{i^N_m-i^r_l \over (i^N_m-i^{N-1}_m)(i^{N-1}_t-i^r_u)}
{\prod_{s\neq l}^m(i^{N-1}_t-i^r_s)
\prod_{j\not\in\La_r,j\neq i^N_m}(i^r_u-j)
\over
\prod_{j\neq i^N_m,i^{N-1}_m,i^{N-1}_t}(i^{N-1}_t-j)}.
\nonumber
\end{eqnarray}
The contribution to the rhs of (\ref{KZeq}) is
\begin{eqnarray}
&&
{(-1)^{u+m}A_r 
\prod^m_{s\neq u}(i^{N-1}_t-i^r_s)
\prod_{j\not\in\La_r,j\neq i^N_m}(i^r_u-j)
\over
(i^N_m-i^{N-1}_m)
\prod_{j\neq i^N_m,i^{N-1}_m,i^{N-1}_t}(i^{N-1}_t-j)}
\sum_{l\neq u}^m
{\prod_{s=1}^{m-1}(i^r_l-i^N_s)
\over
(i^N_m-i^r_l)(i^{N-1}_t-i^r_l)
\prod_{s\neq l}^m(i^r_l-i^r_s)}.
\label{cont51}
\end{eqnarray}
\vskip5mm

\noindent
(V-II)
The contribution from the term which contains
$\Omega(r,l,\La,k)$ with $k\neq l$.\par
If $l\neq u$ then $\det^1_{IJ}=0$.
In fact if further $k\neq u$ then $k$-th, $l$-th and $u$-th 
rows are proportional and if $k=u$ then $l$-th row is a null vector.
Thus we assume $l=u$.
Then
\begin{eqnarray}
\hbox{$\det^1_{IJ}$}&=&
(-1)^{u+m-1}
{i^N_m-i^r_u \over i^N_m-i^{N-1}_m}
{\prod_{s\neq u}^m(i^{N-1}_t-i^r_s)
\over
\prod_{j\neq i^N_m,i^{N-1}_m,i^{N-1}_t}(i^{N-1}_t-j)}
{1\over \prod_{s=1}^mG^r_s}
F^r_{kt}\prod_{s\neq k}^mG^r_s
\nonumber
\\
&=&(-1)^{u+m-1}
{i^N_m-i^r_u
\over (i^N_m-i^{N-1}_m)(i^{N-1}_t-i^r_k)}
{\prod_{s\neq u}^m(i^{N-1}_t-i^r_s)
\prod_{j\not\in\La_r,j\neq i^N_m}(i^r_k-j)
\over
\prod_{j\neq i^N_m,i^{N-1}_m,i^{N-1}_t}(i^{N-1}_t-j)}.
\nonumber
\end{eqnarray}
The contribution to the rhs of (\ref{KZeq}) is
\begin{eqnarray}
&&
{(-1)^{u+m-1}A_r\prod_{s\neq u}^m(i^{N-1}_t-i^r_s)
\prod_{j\not\in\La_r,j\neq i^N_m,i^{N-1}_m}(i^r_u-j)
\over
(i^N_m-i^{N-1}_m)
\prod_{j\neq i^N_m,i^{N-1}_m,i^{N-1}_t}(i^{N-1}_t-j)}
\nonumber
\\
&&
\times
\sum_{k\neq u}^m
{(i^r_k-i^{N-1}_m)
\prod_{s=1}^{m-1}(i^r_k-i^N_s)
\over 
(i^N_m-i^r_k)(i^{N-1}_t-i^r_k)
\prod_{s\neq k}^m(i^r_k-i^r_s)}.
\label{cont52}
\end{eqnarray}
Note that the $k=u$ term of this equation is equal to the minus
of the $l=u$ term in (\ref{cont51}).
\vskip5mm

\noindent
(V-III)
The contribution from the term which contains
$\Omega(r,l,N-1,k)$.\par
If $k\neq t$, $\det^2_{IJ}=0$. We assume $k=t$.
Then $\det^2_{IJ}=0$ for $l\neq u$, since $l$-th and $u$-th rows are 
proportional. Thus we assume $l=u$. 
Then
\begin{eqnarray}
&&
\hbox{$\det^2_{IJ}$}=(-1)^{u+m}.
\nonumber
\end{eqnarray}
The contribution to the rhs of (\ref{KZeq}) is
\begin{eqnarray}
&&
{(-1)^{u+m}A_r\prod_{j\not\in\La_r,j\neq i^N_m}(i^r_u-j)
\over
(i^N_m-i^r_u)(i^r_u-i^{N-1}_m)(i^N_m-i^{N-1}_t)
\prod_{j\in{\cal K},j\neq i^{N-1}_t}(i^{N-1}_t-j)}
B(r,u,t)
\nonumber
\\
&&=
{(-1)^{u+m}A_r\prod_{j\not\in\La_r,j\neq i^N_m}(i^r_u-j)
\over
(i^N_m-i^r_u)(i^N_m-i^{N-1}_m)(i^N_m-i^{N-1}_t)
\prod_{j\in{\cal K},j\neq i^{N-1}_t}(i^{N-1}_t-j)}
B(r,u,t)
\nonumber
\\
&&
+
{(-1)^{u+m}A_r\prod_{j\not\in\La_r,j\neq i^N_m}(i^r_u-j)
\over
(i^r_u-i^{N-1}_m)(i^N_m-i^{N-1}_m)(i^N_m-i^{N-1}_t)
\prod_{j\in{\cal K},j\neq i^{N-1}_t}(i^{N-1}_t-j)}
B(r,u,t),
\label{cont53}
\end{eqnarray}
where we have used
\begin{eqnarray}
&&
{1\over (i^N_m-i^r_u)(i^r_u-i^{N-1}_m)}
=
{1\over (i^N_m-i^r_u)(i^N_m-i^{N-1}_m)}
+
{1\over (i^r_u-i^{N-1}_m)(i^N_m-i^{N-1}_m)}.
\nonumber
\end{eqnarray}
Note that the first term of (\ref{cont53}) is nothing
but the coefficient of $\Omega^r_{IJ}$ in (\ref{calcueq}).
\vskip2mm

Let us calculate $\hbox{(\ref{cont51})}+\hbox{(\ref{cont52})}$.
We have
\begin{eqnarray}
&&
\hbox{(\ref{cont51})}+\hbox{(\ref{cont52})}
\nonumber
\\
&&=
{(-1)^{u+m}A_r 
\prod^m_{s\neq u}(i^{N-1}_t-i^r_s)
\prod_{j\not\in\La_r,j\neq i^N_m}(i^r_u-j)
\over
(i^N_m-i^{N-1}_m)
\prod_{j\neq i^N_m,i^{N-1}_m,i^{N-1}_t}(i^{N-1}_t-j)}
\nonumber
\\
&&
\times\sum_{l=1}^m
\Big[
1+{i^r_l-i^{N-1}_m \over i^r_u-i^{N-1}_m}
\Big]
{\prod_{s=1}^{m-1}(i^r_l-i^N_s)
\over
(i^N_m-i^r_l)(i^{N-1}_t-i^r_l)
\prod_{s\neq l}^m(i^r_l-i^r_s)}.
\nonumber
\end{eqnarray}
By the residue theorem
\begin{eqnarray}
&&
\sum_{l=1}^m
{\prod_{s=1}^{m-1}(i^r_l-i^N_s)
\over
(i^N_m-i^r_l)(i^{N-1}_t-i^r_l)
\prod_{s\neq l}^m(i^r_l-i^r_s)}
\nonumber
\\
&&=
{\prod_{s=1}^{m-1}(i^N_m-i^N_s)
\over
(i^{N-1}_t-i^N_m)\prod_{s=1}^m(i^N_m-i^r_s)}
+
{\prod_{s=1}^{m-1}(i^{N-1}_t-i^N_s)
\over
(i^N_m-i^{N-1}_t)\prod_{s=1}^m(i^{N-1}_t-i^r_s)},
\nonumber
\\
&&
\sum_{k=1}^m
{(i^r_k-i^{N-1}_m)
\prod_{s=1}^{m-1}(i^r_k-i^N_s)
\over 
(i^N_m-i^r_k)(i^{N-1}_t-i^r_k)
\prod_{s\neq k}^m(i^r_k-i^r_s)}
\nonumber
\\
&&=
{(i^N_m-i^{N-1}_m)\prod_{s=1}^{m-1}(i^N_m-i^N_s)
\over
(i^{N-1}_t-i^N_m)\prod_{s=1}^m(i^N_m-i^r_s)}
+
{(i^{N-1}_t-i^{N-1}_m)\prod_{s=1}^{m-1}(i^{N-1}_t-i^N_s)
\over
(i^N_m-i^{N-1}_t)\prod_{s=1}^m(i^{N-1}_t-i^r_s)}.
\nonumber
\end{eqnarray}
Hence
\begin{eqnarray}
&&
\hbox{(\ref{cont51})}+\hbox{(\ref{cont52})}
\nonumber
\\
&&=
{(-1)^{u+m}A_r 
\prod_{s\neq u}(i^{N-1}_t-i^r_s)
\prod_{j\not\in\La_r,j\neq i^N_m}(i^r_u-j)
\over
(i^N_m-i^{N-1}_m)(i^N_m-i^{N-1}_t)
\prod_{j\neq i^N_m,i^{N-1}_m,i^{N-1}_t}(i^{N-1}_t-j)}
\times
\nonumber
\\
&&
\times
\Big[
{i^N_m-i^r_u \over i^r_u-i^{N-1}_m}
{\prod_{s=1}^{m-1}(i^N_m-i^N_s) \over \prod_{s=1}^{m}(i^N_m-i^r_s)}
+{i^r_u-i^{N-1}_t \over i^r_u-i^{N-1}_m}
{\prod_{s=1}^{m-1}(i^{N-1}_t-i^N_s) \over \prod_{s=1}^{m}(i^{N-1}_t-i^r_s)}
\Big]
\nonumber
\\
&&=
{(-1)^{u+m+1}A_r 
\prod_{j\not\in\La_r,j\neq i^N_m}(i^r_u-j)
\over
(i^r_u-i^{N-1}_m)(i^N_m-i^{N-1}_m)(i^N_m-i^{N-1}_t)
\prod_{j\in{\cal K},j\neq i^{N-1}_t}(i^{N-1}_t-j)}
B(r,u,t).
\nonumber
\end{eqnarray}
This is the minus of the second term of (\ref{cont53}).
Hence
\begin{eqnarray}
&&
\hbox{(\ref{cont51})}+\hbox{(\ref{cont52})}+\hbox{(\ref{cont53})}
=\hbox{ the first term of (\ref{cont53})}
\nonumber
\end{eqnarray}
which is equal to the coefficient of $\Omega^r_{IJ}$ in (\ref{calcueq}).
\vskip7mm

\noindent
{\bf(VI)}: The coefficient of 
$\Omega^{N-1}_{IJ}$ for which $I\cap J=\emptyset$,
$|I|\geq 1$ and $|J|\geq 2$.
\vskip3mm
Let us set 
$I=\{p_1<\cdots<p_u\}$ and $J=\{q_1<\cdots<q_t\}$, $u,t\leq m-1$, $u+t=m-1$.
\vskip2mm
\noindent
(VI-I)
The contribution from the term which contains
$\Omega(N-1,l,\La,l)$.\par
If $l\not\in I$ then $\det^3_{IJ}=0$. 
We assume $l\in I$. 
We have
\begin{eqnarray}
\hbox{$\det^3_{IJ}$}&=&
\sgn\cdot \prod_{s\in I\backslash\{l\}}C^l_{ss}
\prod_{s\in J}D^l_{ss}
\nonumber
\\
&=&
\sgn\cdot
{(i^{N-1}_l-i^{N-1}_m)^{u-1}(i^N_m-i^{N-1}_l)^t 
\over (i^N_m-i^{N-1}_m)^{m-2}}
\prod_{s\in I\backslash\{l\}}
{i^N_m-i^{N-1}_s \over i^{N-1}_l-i^{N-1}_s}
\prod_{s\in J}
{i^{N-1}_m-i^{N-1}_s \over i^{N-1}_l-i^{N-1}_s},
\nonumber
\end{eqnarray}
where $\sgn=\sgn(p_1,\cdots,p_u,q_1,\cdots,q_t)$ is the sign
of the permutation.
The contribution to the rhs of (\ref{KZeq}) is
\begin{eqnarray}
&&
\sgn
{\prod_{s\in I}(i^N_m-i^{N-1}_s)
\prod_{s\in J}(i^{N-1}_m-i^{N-1}_s)
\over (i^N_m-i^{N-1}_m)^{m-2}}
\sum_{l\in I}
{(i^{N-1}_l-i^{N-1}_m)^{u-1}(i^N_m-i^{N-1}_l)^{t-2}
\over 
\prod_{s\neq l}^{m-1}(i^{N-1}_l-i^{N-1}_s)}.
\label{cont61}
\end{eqnarray}
\vskip2mm

\noindent
(VI-II)
The contribution from the term which contains
$\Omega(N-1,l,N-1,l)$ 
with $1\leq l\leq m-1$.\par
If $l\not\in J$, $\det^4_{IJ}=0$. 
We assume $l\in J$.
We have
\begin{eqnarray}
\hbox{$\det^4_{IJ}$}&=&
\sgn\cdot
\prod_{s\in I}C^l_{ss}
\prod_{s\in J\backslash\{l\}}D^l_{ss}
\nonumber
\\
&=&
\sgn\cdot
{(i^{N-1}_l-i^{N-1}_m)^{u}(i^N_m-i^{N-1}_l)^{t-1} 
\over (i^N_m-i^{N-1}_m)^{m-2}}
\prod_{s\in I}
{i^N_m-i^{N-1}_s \over i^{N-1}_l-i^{N-1}_s}
\prod_{s\in J\backslash\{l\}}
{i^{N-1}_m-i^{N-1}_s \over i^{N-1}_l-i^{N-1}_s}.
\nonumber
\end{eqnarray}
The contribution to the rhs of (\ref{KZeq}) is
\begin{eqnarray}
&&
-\sgn
{A_{N-1}\prod_{s\in I}(i^N_m-i^{N-1}_s)
\prod_{s\in J}(i^{N-1}_m-i^{N-1}_s)
\over (i^N_m-i^{N-1}_m)^{m-2}}
\times
\nonumber
\\
&&
\times
\sum_{l\in J}
{(i^{N-1}_l-i^{N-1}_m)^{u-2}(i^N_m-i^{N-1}_l)^{t-3} 
\prod_{s=1}^{m-1}(i^{N-1}_l-i^N_s)
\over
\prod_{s\neq l}^{m-1}(i^{N-1}_l-i^{N-1}_s)^2}
\nonumber
\\
&&+
\sgn
{\prod_{s\in I}(i^N_m-i^{N-1}_s)
\prod_{s\in J}(i^{N-1}_m-i^{N-1}_s)
\over (i^N_m-i^{N-1}_m)^{m-2}}
\sum_{l\in J}
{(i^{N-1}_l-i^{N-1}_m)^{u-1}(i^N_m-i^{N-1}_l)^{t-2} 
\over
\prod_{s\neq l}^{m-1}(i^{N-1}_l-i^{N-1}_s)}.
\label{cont62}
\end{eqnarray}
Note that the sum of the second term of this equation and (\ref{cont61})
is zero by the residue theorem and the conditions $u\geq 1$, $t\geq 2$.
\vskip2mm

\noindent
(VI-III)
The contribution from the term which contains
$\Omega(N-1,l,N-1,k)$ 
with $1\leq k,l\leq m-1$, $k\neq l$.\par
If $k\not\in J$ then $\det^4_{IJ}=0$. 
We assume $k\in J$.
\vskip2mm

\noindent
(VI-III-I) $l\in I$ case.\par
We have
\begin{eqnarray}
\hbox{$\det^4_{IJ}$}
&=&
-\sgn \cdot C^l_{kl}
\prod_{s\in I\backslash\{l\}}C^l_{ss}
\prod_{s\in J\backslash\{k\}}D^l_{ss}
\nonumber
\\
&=&
\sgn
{ (i^{N-1}_l-i^{N-1}_m)^{u}(i^N_m-i^{N-1}_l)^{t} 
\prod_{s\in I\backslash\{l\}}(i^N_m-i^{N-1}_s)
\prod_{s\in J\backslash\{k\}}(i^{N-1}_m-i^{N-1}_s)
\over (i^N_m-i^{N-1}_m)^{m-2} \prod_{s\neq l}^{m-1}(i^{N-1}_l-i^{N-1}_s)}
\nonumber
\\
&&
\times
\prod_{j\not\in\La_{N-1},j\neq i^N_m}
{i^{N-1}_k-j \over i^{N-1}_l-j}.
\nonumber
\end{eqnarray}
The contribution to the rhs of (\ref{KZeq}) is
\begin{eqnarray}
&&
\sgn
{ A_{N-1}
\prod_{s\in I}(i^N_m-i^{N-1}_s)
\prod_{s\in J}(i^{N-1}_m-i^{N-1}_s) 
\over (i^N_m-i^{N-1}_m)^{m-2} }
\sum_{l\in I}
{(i^{N-1}_l-i^{N-1}_m)^{u-1}(i^N_m-i^{N-1}_l)^{t-2} 
\over
\prod_{s\neq l}^{m-1}(i^{N-1}_l-i^{N-1}_s)}
\nonumber
\\
&&
\times
\sum_{k\in J}
{\prod_{s=1}^{m-1}(i^{N-1}_k-i^N_s)
\over
(i^{N-1}_k-i^N_m)\prod_{s\neq k}^m(i^{N-1}_k-i^{N-1}_s)}.
\label{cont631}
\end{eqnarray}
\vskip2mm

\noindent
(VI-III-II) $l\in J$ case.\par
We have
\begin{eqnarray}
\hbox{$\det^4_{IJ}$}&=&
-\sgn \cdot  D^l_{kl}
\prod_{s\in I}C^l_{ss}
\prod_{s\in J\backslash\{k,l\}}D^l_{ss}
\nonumber
\\
&=&
-\sgn
{(i^{N-1}_l-i^{N-1}_m)^{u+1}(i^N_m-i^{N-1}_l)^{t-1} 
\prod_{s\in I}(i^N_m-i^{N-1}_s)
\prod_{s\in J\backslash\{k,l\}}(i^{N-1}_m-i^{N-1}_s) 
\over 
(i^N_m-i^{N-1}_m)^{m-2}
\prod_{s\neq l}^{m-1}(i^{N-1}_l-i^{N-1}_s)
}
\nonumber
\\
&&
\times
\prod_{j\not\in\La_{N-1},j\neq i^N_m}
{i^{N-1}_k-j \over i^{N-1}_l-j}.
\nonumber
\end{eqnarray}
The contribution to the rhs of (\ref{KZeq}) is
\begin{eqnarray}
&&
\sgn
{A_{N-1}
\prod_{s\in I}(i^N_m-i^{N-1}_s)
\prod_{s\in J}(i^{N-1}_m-i^{N-1}_s) 
\over (i^N_m-i^{N-1}_m)^{m-2}}
\sum_{l\in J}
{(i^{N-1}_l-i^{N-1}_m)^{u-1}(i^N_m-i^{N-1}_l)^{t-2} 
\over
\prod_{s\neq l}^{m-1}(i^{N-1}_l-i^{N-1}_s)}
\nonumber
\\
&&
\times
\sum_{k\in J,k\neq l}
{\prod_{s=1}^{m-1}(i^{N-1}_k-i^N_s)
\over
(i^{N-1}_k-i^N_m)\prod_{s\neq k}^m(i^{N-1}_k-i^{N-1}_s)}.
\label{cont632}
\end{eqnarray}
If we set $k=l$ in this equation, then it is equal to the
first term of (\ref{cont62}).
\vskip2mm

Thus we have
\begin{eqnarray}
&&
\hbox{(\ref{cont61})}+\hbox{(\ref{cont62})}+\hbox{(\ref{cont631})}
+\hbox{(\ref{cont632})}
\nonumber
\\
&&=
\sgn
{A_{N-1}
\prod_{s\in I}(i^N_m-i^{N-1}_s)
\prod_{s\in J}(i^{N-1}_m-i^{N-1}_s) 
\over (i^N_m-i^{N-1}_m)^{m-2}}
\sum_{l=1}^{m-1}
{(i^{N-1}_l-i^{N-1}_m)^{u-1}(i^N_m-i^{N-1}_l)^{t-2} 
\over
\prod_{s\neq l}^{m-1}(i^{N-1}_l-i^{N-1}_s)}
\nonumber
\\
&&
\times
\sum_{k\in J}
{\prod_{s=1}^{m-1}(i^{N-1}_k-i^N_s)
\over
(i^{N-1}_k-i^N_m)\prod_{s\neq k}^m(i^{N-1}_k-i^{N-1}_s)}
\nonumber
\\
&&=0,
\nonumber
\end{eqnarray}
by applying the residue theorem to the summation in $l$.
\vskip7mm

\noindent
{\bf(VII)}: The coefficient of 
$\Omega^{N-1}_{IJ}$ for which $|I\cap J|=1$,
$|J|\geq 2$.
\vskip3mm
\noindent
We set $I=\{p_1<\cdots<p_u\}$, $J=\{q_1<\cdots<q_t\}$ $(u+t=m-1)$
and $I\cap J=\{\bk\}$.\par

\noindent
(VII-I)
The contribution from the term which contains
$\Omega(N-1,l,\La,l)$.\par
If $\bk\neq l$ then $\det^3_{IJ}=0$. 
In fact if $\bk\neq l$, either the $l$-th row is a null vector or
there exists a row proportional to the $l$-th row.
We assume $l=\bk$. Let us define $v,w,y$ by $p_v=q_w=l$,
$I\cup J=\{1,2,\cdots,m-1\}\backslash\{y\}$.
Then we have
\begin{eqnarray}
\hbox{$\det^3_{IJ}$}&=&
\sgn\cdot D^l_{yl}\prod_{s\in I\backslash\{l\}}C^l_{ss}
\prod_{s\in J\backslash\{l\}}D^l_{ss}
\nonumber
\\
&=&
\sgn\cdot
{(i^{N-1}_l-i^{N-1}_m)^{u}(i^N_m-i^{N-1}_l)^t 
\prod_{s\in I\backslash\{l\}}(i^N_m-i^{N-1}_s)
\prod_{s\in J\backslash\{l\}}(i^{N-1}_m-i^{N-1}_s)
\over (i^N_m-i^{N-1}_m)^{m-2}
\prod_{s\neq l}^{m-1}(i^{N-1}_l-i^{N-1}_s)}
\nonumber
\\
&&
\times
\prod_{j\not\in\La_{N-1},j\neq i^N_m}
{i^{N-1}_y-j \over i^{N-1}_l-j},
\nonumber
\end{eqnarray}
where $\sgn=\sgn(p_1,\cdots,p_u,q_1,\cdots,y,\cdots,q_t)$, $y$ being
on the place of $q_w$.
The contribution to the rhs of (\ref{KZeq}) is
\begin{eqnarray}
&&
-\sgn\cdot
{(i^{N-1}_{\bk}-i^{N-1}_m)^{u-1}(i^N_m-i^{N-1}_{\bk})^{t-2} 
\prod_{s\in I}(i^N_m-i^{N-1}_s)
\prod_{s\in J}(i^{N-1}_m-i^{N-1}_s)
\over (i^N_m-i^{N-1}_m)^{m-2}
\prod_{s\neq \bk}^{m-1}(i^{N-1}_{\bk}-i^{N-1}_s)}
\nonumber
\\
&&
\times
\prod_{j\not\in\La_{N-1},j\neq i^N_m}
{i^{N-1}_y-j \over i^{N-1}_{\bk}-j}.
\label{cont71}
\end{eqnarray}
\vskip2mm

\noindent
(VII-II)
The contribution from the term which contains
$\Omega(N-1,l,N-1,l)$ 
with $1\leq l\leq m-1$.\par
If $\bk\neq l$ then $\det^4_{IJ}=0$ by the same reason as (VII-I). 
We assume $l=\bk$. Let us define $v,w,y$ as in (VII-I).
We have
\begin{eqnarray}
\hbox{$\det^4_{IJ}$}&=&
-\sgn\cdot C^l_{yl}
\prod_{s\in I\backslash\{l\}}C^l_{ss}
\prod_{s\in J\backslash\{l\}}D^l_{ss}
\nonumber
\\
&=&
\sgn\cdot D^l_{yl}
\prod_{s\in I\backslash\{l\}}C^l_{ss}
\prod_{s\in J\backslash\{l\}}D^l_{ss}
\nonumber
\end{eqnarray}
which is same as $\det^3_{IJ}$ in (VII-I).
The contribution to the rhs of (\ref{KZeq}) is
\begin{eqnarray}
&&
-\sgn
{
A_{N-1}
(i^{N-1}_{\bk}-i^{N-1}_m)^{u-2}(i^N_m-i^{N-1}_{\bk})^{t-3} 
\over 
(i^N_m-i^{N-1}_m)^{m-2}
}
\times
\nonumber
\\
&&
\times
{
\prod_{s=1}^{m-1}(i^{N-1}_{\bk}-i^N_s)
\prod_{s\in I}(i^N_m-i^{N-1}_s)
\prod_{s\in J}(i^{N-1}_m-i^{N-1}_s)
\over
\prod_{s\neq \bk}^{m-1}(i^{N-1}_{\bk}-i^{N-1}_s)^2
}
\prod_{j\not\in\La_{N-1},j\neq i^N_m}
{i^{N-1}_y-j \over i^{N-1}_{\bk}-j}
\nonumber
\\
&&
+\sgn\cdot
{(i^{N-1}_{\bk}-i^{N-1}_m)^{u-1}(i^N_m-i^{N-1}_{\bk})^{t-2} 
\prod_{s\in I}(i^N_m-i^{N-1}_s)
\prod_{s\in J}(i^{N-1}_m-i^{N-1}_s)
\over (i^N_m-i^{N-1}_m)^{m-2}
\prod_{s\neq \bk}^{m-1}(i^{N-1}_{\bk}-i^{N-1}_s)}
\times
\nonumber
\\
&&
\times
\prod_{j\not\in\La_{N-1},j\neq i^N_m}
{i^{N-1}_y-j \over i^{N-1}_{\bk}-j}.
\label{cont72}
\end{eqnarray}
The second term of this equation is equal to the minus of
(\ref{cont71}).
\vskip2mm

\noindent
(VII-III)
The contribution from the term which contains
$\Omega(N-1,l,N-1,k)$ 
with $1\leq k,l\leq m-1$, $k\neq l$.\par
If $\bk\neq k$ then $\det^4_{IJ}=0$. 
In fact if $\bk\neq k$ then $\bk$-th column in the right matrix and
that in the left matrix are proportional.
We assume $k=\bk$.
Let us define $v,w,y$ by $p_v=q_w=k$,
$I\cup J=\{1,2,\cdots,m-1\}\backslash\{y\}$.
\vskip2mm

\noindent
(VII-III-I) $l\in I$ case.\par
We have
\begin{eqnarray}
\hbox{$\det^4_{IJ}$}&=&
-\sgn\cdot C^l_{yl}\prod_{s\in I\backslash\{l\}}C^l_{ss}
\prod_{s\in J\backslash\{k\}}D^l_{ss}
\nonumber
\\
&=&
\sgn\cdot
{(i^{N-1}_l-i^{N-1}_m)^{u}(i^N_m-i^{N-1}_l)^t 
\prod_{s\in I\backslash\{l\}}(i^N_m-i^{N-1}_s)
\prod_{s\in J\backslash\{k\}}(i^{N-1}_m-i^{N-1}_s)
\over (i^N_m-i^{N-1}_m)^{m-2}
\prod_{s\neq l}^{m-1}(i^{N-1}_l-i^{N-1}_s)}
\nonumber
\\
&&
\prod_{j\not\in\La_{N-1},j\neq i^N_m}
{i^{N-1}_y-j \over i^{N-1}_l-j}.
\nonumber
\end{eqnarray}
The contribution to the rhs of (\ref{KZeq}) is
\begin{eqnarray}
&&
\sgn \cdot A_{N-1}
{\prod_{s\in I\backslash\{\bk\}}(i^N_m-i^{N-1}_s)
\prod_{s\in J\backslash\{\bk\}}(i^{N-1}_m-i^{N-1}_s)
 \prod_{j\not\in\La_{N-1},j\neq i^N_m}(i^{N-1}_y-j)
\over
(i^N_m-i^{N-1}_m)^{m-2}
\prod_{j\in{\cal K},j\neq i^{N-1}_{\bk}}(i^{N-1}_{\bk}-j)}
\nonumber
\\
&&
\times
\sum_{l\in I,l\neq \bk}
{
(i^{N-1}_l-i^{N-1}_m)^{u-1}(i^N_m-i^{N-1}_l)^{t-2} 
\over
\prod_{s\neq l}^{m-1}(i^{N-1}_l-i^{N-1}_s)}
\label{cont731}
\end{eqnarray}
\vskip2mm

\noindent
(VII-III-II) $l\in J$ case.\par
We have
\begin{eqnarray}
\hbox{$\det^4_{IJ}$}&=&
-\sgn\cdot D^l_{yl}\prod_{s\in I}C^l_{ss}
\prod_{s\in J\backslash\{k,l\}}D^l_{ss}
\nonumber
\\
&=&
\sgn\cdot
{(i^{N-1}_l-i^{N-1}_m)^{u}(i^N_m-i^{N-1}_l)^{t-1}
\prod_{s\in I}(i^N_m-i^{N-1}_s)
\prod_{s\in J\backslash\{k\}}(i^{N-1}_m-i^{N-1}_s)
\over (i^N_m-i^{N-1}_m)^{m-2}
\prod_{s\neq l}^{m-1}(i^{N-1}_l-i^{N-1}_s)}
\nonumber
\\
&&
\times
\prod_{j\not\in\La_{N-1},j\neq i^N_m}
{i^{N-1}_y-j \over i^{N-1}_l-j}.
\nonumber
\end{eqnarray}
The contribution to the rhs of (\ref{KZeq}) is
\begin{eqnarray}
&&
\sgn\cdot A_{N-1}
{\prod_{s\in I\backslash\{\bk\}}(i^N_m-i^{N-1}_s)
\prod_{s\in J\backslash\{\bk\}}(i^{N-1}_m-i^{N-1}_s)
\prod_{j\not\in\La_{N-1},j\neq i^N_m}(i^{N-1}_y-j)
\over
(i^N_m-i^{N-1}_m)^{m-2}
\prod_{j\in{\cal K},j\neq i^{N-1}_{\bk}}(i^{N-1}_{\bk}-j)}
\nonumber
\\
&&
\times
\sum_{l\in J,l\neq \bk}
{
(i^{N-1}_l-i^{N-1}_m)^{u-1}(i^N_m-i^{N-1}_l)^{t-2} 
\over
\prod_{s\neq l}^{m-1}(i^{N-1}_l-i^{N-1}_s)}.
\label{cont732}
\end{eqnarray}
Note that the first term of (\ref{cont72}) is equal to
the $l=\bk$ term of (\ref{cont732}).
\vskip2mm

\noindent
(VII-III-III) $l=y$ case.\par
We have
\begin{eqnarray}
\hbox{$\det^4_{IJ}$}&=&
\sgn\cdot \prod_{s\in I}C^l_{ss}
\prod_{s\in J\backslash\{k\}}D^l_{ss}
\nonumber
\\
&=&
\sgn\cdot
{(i^{N-1}_l-i^{N-1}_m)^{u}(i^N_m-i^{N-1}_l)^{t-1}
\prod_{s\in I}(i^N_m-i^{N-1}_s)
\prod_{s\in J\backslash\{k\}}(i^{N-1}_m-i^{N-1}_s)
\over (i^N_m-i^{N-1}_m)^{m-2}
\prod_{s\neq l}^{m-1}(i^{N-1}_l-i^{N-1}_s)}.
\nonumber
\end{eqnarray}
The contribution to the rhs of (\ref{KZeq}) is
\begin{eqnarray}
&&
\sgn\cdot A_{N-1}
{\prod_{s\in I\backslash\{\bk\}}(i^N_m-i^{N-1}_s)
\prod_{s\in J\backslash\{\bk\}}(i^{N-1}_m-i^{N-1}_s)
\prod_{j\not\in\La_{N-1},j\neq i^N_m}(i^{N-1}_y-j)
\over
(i^N_m-i^{N-1}_m)^{m-2}
\prod_{j\in{\cal K},j\neq i^{N-1}_{\bk}}(i^{N-1}_{\bk}-j)}
\nonumber
\\
&&
\times
{
(i^{N-1}_y-i^{N-1}_m)^{u-1}(i^N_m-i^{N-1}_y)^{t-2} 
\over
\prod_{s\neq y}^{m-1}(i^{N-1}_y-i^{N-1}_s)}.
\label{cont733}
\end{eqnarray}
This equation coinsides with that obtained from (\ref{cont731})
or (\ref{cont732}) by setting $l=y$.
\vskip2mm

Now we have
\begin{eqnarray}
&&
\hbox{(\ref{cont71})}+\hbox{(\ref{cont72})}+
\hbox{(\ref{cont731})}+\hbox{(\ref{cont732})}+
\hbox{(\ref{cont733})}
\nonumber
\\
&&=
\sgn\cdot A_{N-1}
{\prod_{s\in I\backslash\{\bk\}}(i^N_m-i^{N-1}_s)
\prod_{s\in J\backslash\{\bk\}}(i^{N-1}_m-i^{N-1}_s)
\prod_{j\not\in\La_{N-1},j\neq i^N_m}(i^{N-1}_y-j)
\over
(i^N_m-i^{N-1}_m)^{m-2}
\prod_{j\in{\cal K},j\neq i^{N-1}_{\bk}}(i^{N-1}_{\bk}-j)}
\nonumber
\\
&&
\times
\sum_{l=1}^{m-1}
{
(i^{N-1}_l-i^{N-1}_m)^{u-1}(i^N_m-i^{N-1}_l)^{t-2} 
\over
\prod_{s\neq l}^{m-1}(i^{N-1}_l-i^{N-1}_s)}
\nonumber
\\
&&=0.
\nonumber
\end{eqnarray}
The last equality follows from the residue theorem.
\vskip7mm

\noindent
{\bf(VIII)}: The coefficient of 
$\Omega^{N-1}_{IJ}$ for which $|I\cap J|\geq 2$.
\vskip3mm

\noindent
(VIII-I)
The contribution from the term which contains
$\Omega(N-1,l,\La,l)$ is zero.\par
In fact $\det^3_{IJ}=0$ since at least one pair of common column
is linearly dependent.
\vskip2mm

\noindent
(VIII-II)
The contribution from the term which contains
$\Omega(N-1,l,N-1,l)$ with
$1\leq l\leq m-1$ is zero by the same reason as (VIII-I).
\vskip2mm

\noindent
(VIII-III)
The contribution from the term which contains
$\Omega(N-1,l,N-1,k)$ with
$1\leq k,l\leq m-1$, $k\neq l$ is zero by the following reason. 
Since $C^l_{kl}=-D^l_{kl}$ for $k\neq l$, the common column
except $k$-th column is linearly dependent.
Hence $\det^4_{IJ}=0$.
\vskip2mm

As a whole the coefficient of $\Omega^{N-1}_{IJ}$ in the rhs
of (\ref{KZeq}) is zero.
\vskip7mm

\noindent
{\bf(IX)}: The coefficient of 
$\Omega^r_{IJ}$ for which $|J|\geq 2$.
\vskip3mm
\noindent
Let us set $I=\{p_1<\cdots<p_u\}$, $J=\{q_1<\cdots<q_t\}$,
$\bI=\{1,2,\cdots,m\}\backslash I=\{\bp_1<\cdots<\bp_t\}$
with $u+t=m$.
\vskip2mm

\noindent
(IX-I)
The contribution from the term which contains
$\Omega(r,l,\La,l)$.\par
If $l\not\in I$, $\det^1_{IJ}=0$.
We assume $l\in I$.
Then
\begin{eqnarray}
\hbox{$\det^1_{IJ}$}&=&
(-1)^{\sum p_i+{1\over2}u(u+1)}
\prod_{s\in I}G^r_s
\det(F_{\bp_iq_j})
\Big(
{i^N_m-i^r_l \over i^N_m-i^{N-1}_m}
\Big)^t
{E_J \over \prod_{s\in J}(i^{N-1}_s-i^r_l)}
{1\over \prod_{s=1}^mG^r_s}
\nonumber
\\
&=&
\sgn{E_JD_{IJ} \over G^r(\bI)}
\Big(
{i^N_m-i^r_l \over i^N_m-i^{N-1}_m}
\Big)^t
{1\over \prod_{s\in J}(i^{N-1}_s-i^r_l)},
\nonumber
\end{eqnarray}
where $\sgn=(-1)^{\sum p_i+{1\over2}u(u+1)+{1\over2}t(t+1)}$,
$G^r(\bI)=\prod_{s\in \bI}G^r_s$
and
\begin{eqnarray}
&&
D_{IJ}=
(-1)^{{1\over2}t(t+1)}\det(F_{\bp_iq_j})=
{\prod_{\alpha<\beta}(i^r_{\bp_\alpha}-i^r_{\bp_\beta})
\prod_{\alpha<\beta}(i^{N-1}_{q_\alpha}-i^{N-1}_{q_\beta})
\over
\prod_{s\in\bI}\prod_{s^\prime\in J}(i^r_s-i^{N-1}_{s^\prime})},
\nonumber
\\
&&
E_J=
\prod_{y\in J}
\Big(
{\prod_{s=1}^m(i^{N-1}_y-i^r_s)
\over
\prod_{j\neq i^N_m,i^{N-1}_m,i^{N-1}_y}(i^{N-1}_y-j)}
\Big).
\nonumber
\end{eqnarray}
The contribution to the rhs of (\ref{KZeq}) is
\begin{eqnarray}
&&
{\sgn\cdot A_r E_JD_{IJ} \over (i^N_m-i^{N-1}_m)^t G^r(\bI)}
\sum_{l\in I}
{(i^N_m-i^r_l)^{t-2}\prod_{s=1}^{m-1}(i^r_l-i^N_s)
\over
\prod_{s\neq l}^m(i^r_l-i^r_s)
\prod_{s\in J}(i^{N-1}_s-i^r_l)}.
\label{cont91}
\end{eqnarray}
\vskip5mm

\noindent
(IX-II)
The contribution from the term which contains
$\Omega(r,l,\La,k)$ with $k\neq l$.\par
If $l\in I$ or $k\not\in I$ then $\det^1_{IJ}=0$.
In fact if $l\in I$ then $l$-th column in the left matrix is a null vector
and if $k\notin I$ then $l$-th row is a null vector.
We assume $l\not\in I$ and $k\in I$.
Let us set $\tilde{I}=\{\tp_1<\cdots<\tp_t\}=
\bI\backslash \{l\}\sqcup\{k\}$ and $p_i^\ast=\bp_i(\bp_i\neq l)$,
$p_i^\ast=k(\bp_i=l)$. 
Let us define $v,w$ by $p_v=k$, $p_{w-1}<l<p_w$.
Then
\begin{eqnarray}
\hbox{$\det^1_{IJ}$}&=&
(-1)^{\sum p_i+{1\over2}u(u+1)+l-k+v-w}
\Big(
{i^N_m-i^r_l \over i^N_m-i^{N-1}_m}
\Big)^t
{
E_J G^r_l\prod_{s\in I\backslash\{k\}}G^r_s\cdot
\det(F_{\tp_iq_j})
\over 
\prod_{s\in J}(i^{N-1}_s-i^r_l)\prod_{s=1}^mG^r_s
}
\nonumber
\\
&=&(-1)^{\sum p_i+{1\over2}u(u+1)-1}
{E_J\over G^r(\bI)}
\Big(
{i^N_m-i^r_l \over i^N_m-i^{N-1}_m}
\Big)^t
{\det(F_{p^\ast_iq_j}) \over \prod_{s\in J}(i^{N-1}_s-i^r_l)}
{G^r_l \over G^r_k}
\nonumber
\\
&=&
-\sgn {E_JD_{IJ} \over G^r(\bI)}
\Big(
{i^N_m-i^r_l \over i^N_m-i^{N-1}_m}
\Big)^t
{1 \over \prod_{s\in J}(i^{N-1}_s-i^r_k)}
{\prod_{j\not\in\La_r,j\neq i^N_m}(i^r_k-j)
\over
\prod_{j\not\in\La_r,j\neq i^N_m}(i^r_l-j)}
{\prod_{s\in \bI\backslash\{l\}}(i^r_k-i^r_s)
\over
\prod_{s\in \bI\backslash\{l\}}(i^r_l-i^r_s)},
\nonumber
\end{eqnarray}
where we use
\begin{eqnarray}
&&
\det(F_{\tp_iq_j})=(-1)^{k-l+w-v-1}\det(F_{p^\ast_iq_j}).
\nonumber
\end{eqnarray}

The contribution to the rhs of (\ref{KZeq}) is
\begin{eqnarray}
&&
-{\sgn\cdot A_r E_ID_{IJ} \over G^r(\bI)(i^N_m-i^{N-1}_m)^t}
\sum_{l\in \bI}
{(i^N_m-i^r_l)^{t-1} 
\over 
(i^r_l-i^{N-1}_m)\prod_{s\in \bI\backslash\{l\}}(i^r_l-i^r_s)}
\times
\nonumber
\\
&&
\times
\sum_{k\in I}
{(i^r_k-i^{N-1}_m)
\prod_{s=1}^{m-1}(i^r_k-i^N_s)
\over
(i^N_m-i^r_k)
\prod_{s\in (I\sqcup \{l\})\backslash\{k\}}(i^r_k-i^r_s)
\prod_{s\in J}(i^{N-1}_s-i^r_k)}.
\label{cont92}
\end{eqnarray}
If we set $k=l$ in this equation then it equals to the minus of
(\ref{cont91}).
\vskip5mm

\noindent
(IX-III)
The contribution from the term which contains
$\Omega(r,l,N-1,k)$.\par
If $l\in I$ or $k\not\in J$, $\det^2_{IJ}=0$.
In fact if $l\in I$ then $l$-th column in the left matrix is zero and
if $k\notin J$ then $l$-th row is zero.
We assume $l\in\bI$ and $k\in J$.
We define $v,w$ here by $q_v=k$ and $p_w<l<p_{w+1}$.
We note that
\begin{eqnarray}
&&
\sharp\{s|\bp_s<l\}=l-1-w,
\quad
\sharp\{s|q_s<k\}=v-1.
\nonumber
\end{eqnarray}
Using these relations we have
\begin{eqnarray}
\hbox{$\det^2_{IJ}$}
&=&(-1)^{\sum p_i+{1\over2}u(u+1)+l-w+v}
\Big(
{i^N_m-i^r_l \over i^N_m-i^{N-1}_m}
\Big)^{t-1}
{E_J \prod_{j\neq i^N_m,i^{N-1}_m,i^{N-1}_k}^m(i^{N-1}_k-j)
\over 
\prod_{s\neq l}^m(i^{N-1}_k-i^r_s)
\prod_{s\in J}(i^{N-1}_s-i^r_l)
\prod_{s\neq l}^mG^r_s}
\times
\nonumber
\\
&&\quad\times
\prod_{s\in I}G^r_s
\det_{(\bI\backslash\{l\})\times(J\backslash\{k\})}(F_{\bp_iq_j})
\nonumber
\\
&=&
\sgn{E_JD_{IJ} \over G^r(\bI)}
\Big(
{i^N_m-i^r_l \over i^N_m-i^{N-1}_m}
\Big)^{t-1}
{\prod_{j\neq i^N_m,i^{N-1}_m,i^{N-1}_k}^m(i^{N-1}_k-j)
\over
\prod_{j\not\in\La_r,j\neq i^N_m}(i^r_l-j)
\prod_{s\neq l}^m(i^{N-1}_k-i^r_s)}
\times
\nonumber
\\
&&\times
{
\prod_{s\in \bI\backslash\{l\}}(i^r_s-i^{N-1}_k)
\over
\prod_{s\in \bI\backslash\{l\}}(i^r_l-i^r_s)
\prod_{s\in J\backslash\{k\}}(i^{N-1}_k-i^{N-1}_s)
}
\end{eqnarray}
The contribution to the rhs of (\ref{KZeq}) is
\begin{eqnarray}
&&
{(-1)^{t-1}
\sgn\cdot A_r E_JD_{IJ} \over G^r(\bI)(i^N_m-i^{N-1}_m)^{t-1}}
\sum_{l\in\bI}
{
(i^N_m-i^r_l)^{t-2} 
\over 
(i^r_l-i^{N-1}_m)\prod_{s\in\bI\backslash\{l\}}(i^r_l-i^r_s)
}
\times
\nonumber
\\
&&\qquad\times
\sum_{k\in J}
{
\prod_{s=1}^{m-1}(i^{N-1}_k-i^N_s)
\over
(i^N_m-i^{N-1}_k)\prod_{s\in I}(i^{N-1}_k-i^r_s)
\prod_{s\in J\backslash\{k\}}(i^{N-1}_k-i^{N-1}_s)
}
\nonumber
\\
&&-
{\sgn\cdot E_JD_{IJ} \over G^r(\bI)(i^N_m-i^{N-1}_m)^{t-1} }
\sum_{l\in\bI}
{
(i^N_m-i^r_l)^{t-1} 
\over 
(i^r_l-i^{N-1}_m)\prod_{s\in\bI\backslash\{l\}}(i^r_l-i^r_s)
}
\sum_{k\in J}
{
\prod_{s\in\bI\backslash\{l\}}(i^r_s-i^{N-1}_k)
\over
(i^N_m-i^{N-1}_k)\prod_{s\in J\backslash\{k\}}(i^{N-1}_k-i^{N-1}_s)
}
\nonumber
\\
&&=
{
(-1)^t\sgn\cdot A_r E_JD_{IJ}
\over 
G^r(\bI)(i^N_m-i^{N-1}_m)\prod_{s\in\bI}(i^{N-1}_m-i^r_s)
}
\sum_{k\in J}
{
\prod_{s=1}^{m-1}(i^{N-1}_k-i^N_s)
\over
(i^N_m-i^{N-1}_k)\prod_{s\in I}(i^{N-1}_k-i^r_s)
\prod_{s\in J\backslash\{k\}}(i^{N-1}_k-i^{N-1}_s)
}
\nonumber
\\
&&
-
{
\sgn\cdot E_JD_{IJ} \prod_{s\in\bI}(i^r_s-i^N_m)
\over
G^r(\bI)(i^N_m-i^{N-1}_m)
\prod_{s\in J}(i^N_m-i^{N-1}_s)
\prod_{s\in\bI}(i^{N-1}_m-i^r_s)}.
\label{cont93}
\end{eqnarray}
In deriving the last equation we use
\begin{eqnarray}
\sum_{l\in\bI}
{
(i^N_m-i^r_l)^{t-2} 
\over 
(i^r_l-i^{N-1}_m)\prod_{s\in\bI\backslash\{l\}}(i^r_l-i^r_s)
}
&=&
-{(i^N_m-i^{N-1}_m)^{t-2} \over \prod_{s\in\bI}(i^{N-1}_m-i^r_s)},
\label{resth2}
\\
\sum_{k\in J}
{
\prod_{s\in\bI\backslash\{l\}}(i^r_s-i^{N-1}_k)
\over
(i^N_m-i^{N-1}_k)\prod_{s\in J\backslash\{k\}}(i^{N-1}_k-i^{N-1}_s)
}
&=&
{\prod_{s\in\bI\backslash\{l\}}(i^r_s-i^N_m)
\over
\prod_{s\in J}(i^N_m-i^{N-1}_s)}.
\nonumber
\end{eqnarray}
Let us name the first and the second term of (\ref{cont93}) 
$Z$ and $W$ respectively.
\vskip5mm

Now we have
\begin{eqnarray}
&&
\hbox{(\ref{cont91})}+\hbox{(\ref{cont92})}+\hbox{(\ref{cont93})}
\nonumber
\\
&&=
{ \sgn\cdot A_r E_JD_{IJ} \over (i^N_m-i^{N-1}_m)^t G^r(\bI) }
\sum_{l=1}^m
{
(i^N_m-i^r_l)^{t-2}\prod_{s=1}^{m-1}(i^r_l-i^N_s)
\over
\prod_{s\neq l}^m(i^r_l-i^r_s)
\prod_{s\in J}(i^{N-1}_s-i^r_l)
}
\nonumber
\\
&&
-
{\sgn\cdot A_r E_ID_{IJ} \over G^r(\bI)(i^N_m-i^{N-1}_m)^t}
\sum_{l\in \bI}
{(i^N_m-i^r_l)^{t-1} 
\over 
(i^r_l-i^{N-1}_m)\prod_{s\in \bI\backslash\{l\}}(i^r_l-i^r_s)}
\times
\nonumber
\\
&&
\times
\sum_{k\in I\sqcup\{l\}}
{(i^r_k-i^{N-1}_m)
\prod_{s=1}^{m-1}(i^r_k-i^N_s)
\over
(i^N_m-i^r_k)
\prod_{s\in (I\sqcup \{l\})\backslash\{k\}}(i^r_k-i^r_s)
\prod_{s\in J}(i^{N-1}_s-i^r_k)}
+Z+W.
\nonumber
\end{eqnarray}
Let us name the fisrt and the second term of this equation
$X$ and $Y$ respectively.
We shall rewrite $X$ and $Y$.
Using
\begin{eqnarray}
&&
\sum_{l=1}^m
{
(i^N_m-i^r_l)^{t-2}\prod_{s=1}^{m-1}(i^r_l-i^N_s)
\over
\prod_{s\neq l}^m(i^r_l-i^r_s)
\prod_{s\in J}(i^{N-1}_s-i^r_l)
}
=
\sum_{l\in J}
{
(i^N_m-i^{N-1}_l)^{t-2}\prod_{s=1}^{m-1}(i^{N-1}_l-i^N_s)
\over
\prod_{s=1}^m(i^{N-1}_l-i^r_s)
\prod_{s\in J\backslash\{l\}}(i^{N-1}_s-i^{N-1}_l)
}
\nonumber
\end{eqnarray}
we have
\begin{eqnarray}
&&
X=
{ \sgn\cdot A_r E_JD_{IJ} \over (i^N_m-i^{N-1}_m)^t G^r(\bI) }
\sum_{l\in J}
{
(i^N_m-i^{N-1}_l)^{t-2}\prod_{s=1}^{m-1}(i^{N-1}_l-i^N_s)
\over
\prod_{s=1}^m(i^{N-1}_l-i^r_s)
\prod_{s\in J\backslash\{l\}}(i^{N-1}_s-i^{N-1}_l)
}.
\nonumber
\end{eqnarray}
Using
\begin{eqnarray}
&&
\sum_{k\in I\sqcup\{l\}}
{(i^r_k-i^{N-1}_m)
\prod_{s=1}^{m-1}(i^r_k-i^N_s)
\over
(i^N_m-i^r_k)
\prod_{s\in (I\sqcup\{l\})\backslash\{k\}}(i^r_k-i^r_s)
\prod_{s\in J}(i^{N-1}_s-i^r_k)}
\nonumber
\\
&&=
{(i^N_m-i^{N-1}_m)\prod_{s=1}^{m-1}(i^N_m-i^N_s)
\over
(i^N_m-i^r_l)\prod_{s\in I}(i^N_m-i^r_s)
\prod_{s\in J}(i^{N-1}_s-i^N_m)}
\nonumber
\\
&&
+
\sum_{k\in J}
{(i^{N-1}_k-i^{N-1}_m)
\prod_{s=1}^{m-1}(i^{N-1}_k-i^N_s)
\over
(i^N_m-i^{N-1}_k)(i^{N-1}_k-i^r_l)
\prod_{s\in I}(i^{N-1}_k-i^r_s)
\prod_{s\in J\backslash\{k\}}(i^{N-1}_s-i^{N-1}_k)}
\nonumber
\end{eqnarray}
we have
\begin{eqnarray}
&&
Y=
-
{\sgn\cdot A_r E_JD_{IJ} 
\prod_{s=1}^{m-1}(i^N_m-i^N_s)
\over G^r(\bI)(i^N_m-i^{N-1}_m)^{t-1}
\prod_{s\in I}^m(i^N_m-i^r_s)
\prod_{s\in J}(i^{N-1}_s-i^N_m)}
\sum_{l\in\bI}
{(i^N_m-i^r_l)^{t-2}
\over
(i^r_l-i^{N-1}_m)
\prod_{s\in\bI\backslash\{l\}}(i^r_l-i^r_s)}
\nonumber
\\
&&
-
{\sgn\cdot A_r E_JD_{IJ} 
\over G^r(\bI)(i^N_m-i^{N-1}_m)^{t}}
\sum_{l\in\bI}
{(i^N_m-i^r_l)^{t-1}
\over
(i^r_l-i^{N-1}_m)
(i^{N-1}_k-i^r_l)
\prod_{s\in\bI\backslash\{l\}}(i^r_l-i^r_s)}
\nonumber
\\
&&
\times
\sum_{k\in J}
{(i^{N-1}_k-i^{N-1}_m)
\prod_{s=1}^{m-1}(i^{N-1}_k-i^N_s)
\over
(i^N_m-i^{N-1}_k)
\prod_{s\in I}(i^{N-1}_k-i^r_s)
\prod_{s\in J\backslash\{k\}}(i^{N-1}_s-i^{N-1}_k)}.
\nonumber
\end{eqnarray}
Using further (\ref{resth2}) and
\begin{eqnarray}
&&
\sum_{l\in\bI}
{(i^N_m-i^r_l)^{t-1}
\over
(i^r_l-i^{N-1}_m)
(i^{N-1}_k-i^r_l)
\prod_{s\in\bI\backslash\{l\}}(i^r_l-i^r_s)}
\nonumber
\\
&=&
-
{(i^N_m-i^{N-1}_m)^{t-1}
\over
(i^{N-1}_k-i^{N-1}_m)\prod_{s\in\bI}(i^{N-1}_m-i^r_s)}
+
{(i^N_m-i^{N-1}_k)^{t-1}
\over
(i^{N-1}_k-i^{N-1}_m)\prod_{s\in\bI}(i^{N-1}_k-i^r_s)}
\nonumber
\end{eqnarray}
we have $Y=-W-Z-X$.
Thus
\begin{eqnarray}
&&
\hbox{(\ref{cont91})}+\hbox{(\ref{cont92})}+\hbox{(\ref{cont93})}
=0.
\end{eqnarray}
This completes the proof of (1) of Theorem \ref{intformula}.

\subsection{Proof of (2) of Theorem 1}
Let us prove the remaining part of Theorem \ref{intformula}.
Let $\{E_{ij}\}$ be the standard basis of $gl_N$ where $E_{ij}$
is the matrix unit with $1$ in $ij$ component.
Set $h_i=E_{ii}-E_{i+1i+1}$ $(1\leq i\leq N-1)$.
By the definition of $f=\sum f_\La v_\La$ it has weight zero,
$h_if=0$ for any $i$.
Hence it is sufficient to prove $E_{ij}f=0$ for any $i\neq j$.
First we assume $N\geq3$.
Then it is sufficient to prove
\begin{eqnarray}
&&
E_{rN}f=0,
\quad
r=1,\cdots, N-2.
\label{singlet1}
\end{eqnarray}
In fact by the following reason the proof for an arbitrary $E_{ij}$ case
is reduced to the above case.
In our description of our basis of differential forms
the index $N$ and $N-1$ play a special role.
For $E_{ij}$ we replace the role of $N$ by $j$ and that of $N-1$ by
$k$ with $k\neq i,j$. This is possible because $N\geq 3$.
Then the following proof is totally the same in this modified
situation.
Thus let us prove (\ref{singlet1}).

Let $\La=(\La_1,\cdots,\La_N)$ with $\La_j=(i^j_1,\cdots,i^j_m)$
$(1\leq j\leq m)$.
We consider $\La^\prime=(\La^\prime_1,\cdots,\La^\prime_N)$ with
\begin{eqnarray}
&&
\La^\prime_j=\La_j \:(j\neq r,N),
\quad
\La^\prime_r=(i^r_1,\cdots,i^r_m,i^N_m),
\quad
\La^\prime_N=(i^N_1,\cdots,i^N_{m-1}).
\nonumber
\end{eqnarray}
Define $v_{\La^\prime}$ in an obvious way.
Then the coefficient of $v_{\La^\prime}$ of $E_{rN}f$ is
\begin{eqnarray}
&&
f_\La+\sum_{l=1}^mf_{\La^{(i^r_li^N_m)}}.
\nonumber
\end{eqnarray}
Hence it is sufficient to prove
\begin{eqnarray}
&&
\bar{f}_\La+\sum_{l=1}^m\bar{f}_{\La^{(i^r_li^N_m)}}=0.
\label{singlet2}
\end{eqnarray}
We shall devide the case into two for the proof of (\ref{singlet2}).
\vskip3mm

\noindent
(i): The coefficient of 
$\Omega=\mdf{\La}{i^1_1}{\La}{i^{N-1}_{m-1}}$ of the left hand side
of (\ref{singlet2}).
\par
\noindent
In (II) of the proof of (1) of Theorem \ref{intformula}
we have calculated the coefficient of $\Omega$ in
\begin{eqnarray}
&&
{\Delta(i^1_1,\cdots,i^{N-1}_{m-1})^{-1} \over 
\prod_{t<u}(\La_t\La_u)}
\sum_{l=1}^m{1\over i^N_m-i^r_l}
\bar{f}_{\La^{(i^r_li^N_m)}}.
\nonumber
\end{eqnarray}
From the calculation there we can easily read off the coefficient
of $\Omega$ in 
$$
{\Delta(i^1_1,\cdots,i^{N-1}_{m-1})^{-1} 
\over \prod_{t<u}(\La_t\La_u)}
\sum_{l=1}^m\bar{f}_{\La^{(i^r_li^N_m)}}.
$$
It is
\begin{eqnarray}
&&
-\sum_{l=1}^m
{A_r \over i^N_m-i^r_l}
{\prod_{s=1}^{m-1}(i^r_l-i^N_s)
\over
\prod_{s\neq l}^m(i^r_l-i^r_s)}
=
-A_r
{\prod_{s=1}^{m-1}(i^N_m-i^N_s)
\over
\prod_{s=1}^m(i^N_m-i^r_s)}
=-1.
\nonumber
\end{eqnarray}
Thus the coefficient of $\Omega$ in the lhs of (\ref{singlet2}) is $0$,
since
$$
{\Delta(i^1_1,\cdots,i^{N-1}_{m-1})^{-1} 
\over \prod_{t<u}(\La_t\La_u)}\bar{f}_\La=\Omega.
$$

\vskip3mm

\noindent
(ii): The coefficient of $\Omega^r_{IJ}$ with $|J|\geq 1$
in the lhs of (\ref{singlet2}).
\par
\noindent
In (IX) we have proved that the coefficient of $\Omega^r_{IJ}$
in
\begin{eqnarray}
&&
{\Delta(i^1_1,\cdots,i^{N-1}_{m-1})^{-1}
\over \prod_{t<u}(\La_t\La_u)}
\sum_{l=1}^m{1\over i^N_m-i^r_l}
\bar{f}_{\La^{(i^r_li^N_m)}}
\nonumber
\end{eqnarray}
is zero. There the condition $|J|\geq 2$ is used only when the residue
theorem is applied. Taking care of it we can again easily read off the
coefficient of $\Omega^r_{IJ}$ 
of
$$
{\Delta(i^1_1,\cdots,i^{N-1}_{m-1})^{-1}
\over \prod_{t<u}(\La_t\La_u)}
\sum_{l=1}^m\bar{f}_{\La^{(i^r_li^N_m)}}
$$
from the calculation in (IX).
We used the notation (\ref{cont91}), (\ref{cont92}), (\ref{cont93}),
$X$, $Y$ to denote the equation apeared there.
We shall use the prime of them for the corresponding equation
like $\hbox{(\ref{cont91})}^\prime$, $X^\prime$ etc.
Then 
\begin{eqnarray}
Z^\prime
&=&
{
(-1)^t\sgn\cdot A_r E_JD_{IJ}
\over 
G^r(\bI)\prod_{s\in\bI}(i^{N-1}_m-i^r_s)
}
\sum_{k\in J}
{
\prod_{s=1}^{m-1}(i^{N-1}_k-i^N_s)
\over
(i^N_m-i^{N-1}_k)\prod_{s\in I}(i^{N-1}_k-i^r_s)
\prod_{s\in J\backslash\{k\}}(i^{N-1}_k-i^{N-1}_s)
}
\nonumber
\\
W^\prime&=&
-{
\sgn\cdot E_JD_{IJ} \prod_{s\in\bI}(i^r_s-i^N_m)
\over
G^r(\bI)
\prod_{s\in J}(i^N_m-i^{N-1}_s)
\prod_{s\in\bI}(i^{N-1}_m-i^r_s)}.
\nonumber
\end{eqnarray}
Also we have
\begin{eqnarray}
&&
X^\prime=
{ \sgn\cdot A_r E_JD_{IJ} \over (i^N_m-i^{N-1}_m)^t G^r(\bI) }
\sum_{l\in J}
{
(i^N_m-i^{N-1}_l)^{t-1}\prod_{s=1}^{m-1}(i^{N-1}_l-i^N_s)
\over
\prod_{s=1}^m(i^{N-1}_l-i^r_s)
\prod_{s\in J\backslash\{l\}}(i^{N-1}_s-i^{N-1}_l)
}
\nonumber
\end{eqnarray}
and
\begin{eqnarray}
&&
Y^\prime=
{
\sgn\cdot E_JD_{IJ} \prod_{s\in\bI}(i^r_s-i^N_m)
\over
G^r(\bI)
\prod_{s\in J}(i^N_m-i^{N-1}_s)
\prod_{s\in\bI}(i^{N-1}_m-i^r_s)}
\nonumber
\\
&&
+
{
(-1)^{t-1}\sgn\cdot A_r E_JD_{IJ}
\over 
G^r(\bI)\prod_{s\in\bI}(i^{N-1}_m-i^r_s)
}
\sum_{k\in J}
{
\prod_{s=1}^{m-1}(i^{N-1}_k-i^N_s)
\over
(i^N_m-i^{N-1}_k)\prod_{s\in I}(i^{N-1}_k-i^r_s)
\prod_{s\in J\backslash\{k\}}(i^{N-1}_k-i^{N-1}_s)
}
\nonumber
\\
&&
-
{ \sgn\cdot A_r E_JD_{IJ} \over (i^N_m-i^{N-1}_m)^t G^r(\bI) }
\sum_{k\in J}
{
(i^N_m-i^{N-1}_k)^{t-1}\prod_{s=1}^{m-1}(i^{N-1}_k-i^N_s)
\over
\prod_{s=1}^m(i^{N-1}_k-i^r_s)
\prod_{s\in J\backslash\{k\}}(i^{N-1}_s-i^{N-1}_k)
}
\nonumber
\\
&&=-W^\prime-Z^\prime-X^\prime.
\end{eqnarray}
Hence
\begin{eqnarray}
&&
\hbox{(\ref{cont91})}^\prime+
\hbox{(\ref{cont92})}^\prime+
\hbox{(\ref{cont93})}^\prime
=Z^\prime+W^\prime+X^\prime+Y^\prime=0.
\nonumber
\end{eqnarray}
Thus the equation (\ref{singlet2}) is proved.

In the $N=2$ case we can similarly read off easily the coefficient
of
$$
{\Delta(i^1_1,\cdots,i^{N-1}_{m-1})^{-1}
\over \prod_{t<u}(\La_t\La_u)}
\big(\bar{f}_\La+\sum_{l=1}^m\bar{f}_{\La^{(i^{N-1}_li^N_m)}}\big)
$$
from (I), (II), (III), (IV), (VI), (VII), (VIII) in the proof of (1) of
Theorem \ref{intformula}.
They are all zero as we expect.
\vskip1cm

\section{Discussion}
In this paper we give integral and theta formulae
for the solutions of $sl_N$ Knizhnik-Zamolodchikov (KZ)
equations of level $0$ with the value in the trivial representation in
the tensor product of the vector representations of $sl_N$. 
The formula generalizes the
Smirnov's formula in the case of $sl_2$.
We have found that the differential form $\mu^\La_p$,
which is a building block of the integral formula,
is obtained by evaluating one of the variables to
the branch point $Q_p$ in the product of chiral
Szeg\"{o} kernels. This is a key for the proof of
the theta formula.

Let us discuss remaining problems and related subjects.

In $N=2$ case it is conjectured that Smirnov type solutions span
the singlet solution space \cite{S4}.
On the other hand the dimension of the vector space spanned by 
our integral formulae is less than the multiplicity of 
the trivial representation in $V^{\ot Nm}$ for $N\geq 3$ and $m\geq2$. 
In fact the multiplicity is given by
$$
\hbox{mult}(0,V^{\ot Nm})={(Nm)! \over
\prod_{k=0}^{N-1}\prod_{j=0}^{m-1}(m+k-j)}.
$$
On the other hand the dimension $D(N,m)$ of the vector space spanned by
integral formulae satisfies
$$
D(N,m)\leq \hbox{I}(N,m)=
\left(\begin{array}{c} Nm-2\\(N-1)m-1\end{array}\right),
$$
where the right hand side is the binomial coefficient.
The number $Nm-2$ is the dimension of an eigenspace of the 
$N$-cyclic automorphism $\phi$ on the first homology group of a $Z_N$ curve.
Then
\begin{eqnarray}
&&
{\hbox{mult}(0,V^{\ot Nm})\over \hbox{I}(N,m)}
=
{N-{1\over m}\over N-1}
\prod_{k=1}^{N-1}\prod_{j=0}^{m-2}
{(N-1)(m-1-j)+k \over m-j+k}.
\label{dratio}
\end{eqnarray}
Since 
$$
(N-1)(m-1-j)+k-(m-j+k)=(N-2)(m-2-j+{N-3\over N-2}),
$$
(\ref{dratio}) is greater than $1$ if $N\geq 3$ and $m\geq 2$.
Note that $\hbox{mult}(0,V^{\ot N})=\hbox{I}(N,1)=1$.
For $N=2$ we have 
$$
\hbox{mult}(0,V^{\ot 2m})=\hbox{I}(2,m)-
\left(\begin{array}{c} 2m-2\\m-3\end{array}\right),
$$
where the second term in the right hand side comes from
the Riemann's bilinear identity \cite{S4}.

This structure of solution space should be same in the qKZ case. 
To construct remaining solutions for both KZ and qKZ equations 
is an interesting and important problem.
In the qKZ case to study a relation 
of these missing solutions with form factors is also interesting.

We still do not understand the relation between the integral formula 
given here and those given in \cite{M,SV1,SV2} 
in the case of $sl_N$, $N\geq 3$. In $N=2$ case the relation 
is given in \cite{NPT}. If we understand this structure then
it will help to find the missing solution discussed above. 

The relation of the solution to the KZ equation of level $0$
with a classical integrable system is still to be clarified.
The relation with the Szeg\"{o} kernel will give some hint
to understand this problem since the Szeg\"{o} kernel
is related with the tau function of the KP hierarchy.
Anyway it is true that we can introduce a Jacobian variable
in the theta formula for the solutions to the KZ equation.
Hence it is natural to ask what kind of equation governs
the dependence on the Jacobian variables and what the zero value
means for that equation.

Once we introduce the Jacobian variable we can ask what is
the difference analogue, $q$ analogue of the theta function?
As to the abelian integral, Smirnov \cite{S4,S5} discussed 
its difference analogue.

Since the Smirnov type formula is related with the 
algebraic curves in the case of $sl_N$, 
it is interesting to study Smirnov type solutions 
for other type of Lie algebras and 
whether they are related with algebraic curves.

The determinantal structure of Smirnov type solution
is still lacking an understanding 
from the representation theoretical view point.

\vskip2cm

\noindent
{\Large\hskip4mm Acknowledgement }
\vskip4mm
\noindent
We would like to thank Koji Cho, Fedor Smirnov, Yasuhiko Yamada
for the stimulating discussions and useful comments.
We also benefited from the discussion with Masaki Kashiwara.
\vskip1cm
\appendix

\section{Appendix 1}
In this section we give a derivation of the formula 
(\ref{rel1})-(\ref{rel5}). We recall the definition of
$\df{\La}{p}$
\begin{eqnarray}
&&
\df{\La}{p}=
{g^{(\La_r)}(p)g^{(p)}_{\La_r}(z)dz \over (z-p)s}
\quad
p\in\La_r,
\nonumber
\end{eqnarray}
where as in the main text $z-p$ means $z-\la_p$ and
$g^{(\La_r)}(p)$ means $g^{(\La_r)}(\la_p)$ etc.

\vskip5mm

\noindent
{\bf (I)}. a derivation of the formula (\ref{rel1}) :
\par
\noindent
By differentiating the defining formula of $\df{\La}{p}$
we have
\begin{eqnarray}
{\partial\over\partial\la_{i^N_m}}
\df{\La}{p}&=&
-{1 \over p-i^N_m}\df{\La}{p}
+{1\over N}
{g^{(\La_r)}(p)g^{(p)}_{\La_r}(z)dz \over (z-p)(z-i^N_m)s}
\label{f11}
\\
&=&
(1-{1\over N}){1 \over i^N_m-p}\df{\La}{p}
+{1\over N}{1 \over i^N_m-p}
{g^{(\La_r)}(p)g^{(p)}_{\La_r}(z)dz \over (z-i^N_m)s}
\label{f12}
\\
&=&
(1-{1\over N}){1 \over i^N_m-p}\df{\La}{p}
-{1\over N}{1 \over i^N_m-p}
\prod_{j\not\in\La_r,j\neq i^N_m}{p-j \over i^N_m-j}
\df{p}{i^N_m}.
\label{f13}
\end{eqnarray}
Here to obtain (\ref{f12}) from (\ref{f11}) we use
\begin{eqnarray}
&&
{1 \over (z-p)(z-i^N_m)}
={1\over i^N_m-p}
\Big[
-{1\over z-p}+{1\over z-i^N_m}
\Big],
\nonumber
\end{eqnarray}
and to get (\ref{f13}) from (\ref{f12}) we use
\begin{eqnarray}
g^{(p)}_{\La_r}(z)&=&
g^{(i^N_m)}_{\La_r^{(pi^N_m)}}(z),
\nonumber
\\
{g^{(\La_r)}(p)
\over
g^{(\La_r^{(pi^N_m)})}(i^N_m)}
&=&
-\prod_{j\not\in\La_r,j\neq i^N_m}{p-j \over i^N_m-j}.
\nonumber
\end{eqnarray}
\vskip5mm

\noindent
{\bf (II)}. a derivation of the formula (\ref{rel2}) :
\par
\noindent
Since
\begin{eqnarray}
&&
\df{i^r_{l'}}{i^r_l}=
{g^{(\La_r^{(i^N_mi^r_{l'})})}(i^r_l)g^{(i^r_l)}_{\La_r^{(i^N_mi^r_{l'})}}(z)dz \over (z-i^r_l)s}
\nonumber
\end{eqnarray}
we have
\begin{eqnarray}
\df{\La}{i^r_l}-\df{i^r_{l'}}{i^r_l}&=&
{g^{(\La_r)}(i^r_l)g^{(i^r_li^r_{l'})}_{\La_r}(z) \over (z-i^r_l)s}dz
\Big[
z-i^r_{l'}
-{i^r_l-i^r_{l'} \over i^r_l-i^N_m}(z-i^N_m)
\Big]
\nonumber
\\
&=&
{i^N_m-i^r_{l'} \over i^N_m-i^r_l}
{g^{(\La_r)}(i^r_l)g^{(i^r_li^r_{l'})}_{\La_r}(z) \over s}dz,
\nonumber
\end{eqnarray}
where we set
$$
g^{(i^r_li^r_{l'})}_{\La_r}(z)=\prod_{j\in\La_r,j\neq i^r_l,i^r_{l'}}(z-j).
$$
Hence
\begin{eqnarray}
&&
{g^{(i^r_li^r_{l'})}_{\La_r}(z) \over s}dz
=
{i^N_m-i^r_{l} \over i^N_m-i^r_{l'}}
{1\over g^{(\La_r)}(i^r_l)}
\Big[
\df{\La}{i^r_l}-\df{i^r_{l'}}{i^r_l}
\Big].
\label{f21}
\end{eqnarray}
If we set $r=N-1$ and $l'=m$ in this equation we get
\begin{eqnarray}
&&
{g^{(i^{N-1}_li^{N-1}_{m})}_{\La_{N-1}}(z) \over s}dz
=
{i^N_m-i^{N-1}_{l} \over i^N_m-i^{N-1}_{m}}
{1\over g^{(\La_{N-1})}(i^{N-1}_l)}
\Big[
\df{\La}{i^{N-1}_l}-\df{i^{N-1}_{m}}{i^{N-1}_l}
\Big].
\label{f22}
\end{eqnarray}
Now let $F(z)$ be an arbitrary polynomial of degree at most $m-2$
then
\begin{eqnarray}
&&
F(z)=\sum_{k=1}^{m-1}
{F(i^{N-1}_k) \over \prod_{s\neq k}^{m-1}(i^{N-1}_k-i^{N-1}_s)}
g^{(i^{N-1}_ki^{N-1}_m)}_{\La_{N-1}}(z).
\label{f23}
\end{eqnarray}
Thus we have
\begin{eqnarray}
&&
{g^{(i^r_li^r_{l'})}_{\La_r}(z) \over s}dz=
\sum_{k=1}^{m-1}
{g^{(i^r_li^r_{l'})}_{\La_r}(i^{N-1}_k) 
\over \prod_{s\neq k}^{m-1}(i^{N-1}_k-i^{N-1}_s)}
{g^{(i^{N-1}_ki^{N-1}_m)}_{\La_{N-1}}(z) \over s}dz.
\label{f24}
\end{eqnarray}
Substitutuing (\ref{f21}) and (\ref{f22}) into (\ref{f24})
we get (\ref{rel2}).
The relation (\ref{rel3}) is a special case $r=N-1$ of (\ref{rel2}).
\vskip5mm

\noindent
{\bf (III)}. Here we prove the formula
\begin{eqnarray}
&&
\df{\La}{i^N_m}=
\sum_{i^r_l\in{\cal K}}
\prod_{j\in{\cal K},j\neq i^r_l}
{i^N_m-j \over i^r_l-j}
\df{\La}{i^r_l}.
\label{f30}
\end{eqnarray}
\par

\noindent
For the sake of simplicity we set $i^r_j=i_{(r-1)m+j}$.
By Proposition \ref{exact} we have
\begin{eqnarray}
&&
[\df{\La}{i_1},\cdots,\df{\La}{i_L}]
=[\zeta^{(\La)}_1,\cdots,\zeta^{(\La)}_L]A,
\nonumber
\end{eqnarray}
where $A$ is the $L$ by $L$ matrix whose $kl$ component $A_{kl}$ is
given by $A_{kl}=\la^{k-1}_{i_l}$.
Let $d_{kl}$ be the $kl$ cofactor of $A$.
Then we have
\begin{eqnarray}
&&
\zeta^{(\La)}_{k}={1\over \det A}\sum_{j=1}^L\df{\La}{i_j}d_{kj}.
\nonumber
\end{eqnarray}
Again by Proposition \ref{exact}
\begin{eqnarray}
&&
\df{\La}{i^N_m}=
\sum_{k=1}^L\zeta^{(\La)}_k\la_{i^N_m}^{k-1}
=
{1\over \det A}
\sum_{j=1}^L\sum_{k=1}^L\la_{i^N_m}^{k-1}d_{kj}
\df{\La}{i_j}.
\label{f31}
\end{eqnarray}
Using the expansion of the Vandermond determinant in a column we have
\begin{eqnarray}
&&
{1\over \det A}
\sum_{k=1}^L\la_{i^N_m}^{k-1}d_{kj}=
\prod_{s\neq j}^L{i^N_m-i_s \over i_j-i_s}.
\label{f32}
\end{eqnarray}
Substitutibg (\ref{f32}) into (\ref{f31}) we have (\ref{f30}).

The relation (\ref{rel5}) is proved in an exactly similar
manner. 
\vskip5mm

\noindent
{\bf (IX)}. a derivation of the formula (\ref{rel4}) :
\par
\noindent
The derivation is similar to (II).
We have
\begin{eqnarray}
&&
\df{\La}{i^N_m}-\df{i^r_l}{i^N_m}
=
g^{(\La_N)}(i^N_m)
{dz \over (z-i^N_m)s}
\Big[
g^{(i^N_m)}_{\La_N}(z)-
{g^{(i^N_m)}_{\La_N}(i^N_m)\over g^{(i^r_l)}_{\La_r}(i^N_m)}
g^{(i^r_l)}_{\La_r}(z)
\Big].
\label{f41}
\end{eqnarray}
Since the polynomial in $[\quad]$ is devided by $z-i^N_m$ we can define
the polynomial $G^{rl}(z)$ by
\begin{eqnarray}
&&
G^{rl}(z)=
{1 \over z-i^N_m}
\Big[
g^{(i^N_m)}_{\La_N}(z)-
{g^{(i^N_m)}_{\La_N}(i^N_m) \over g^{(i^r_l)}_{\La_r}(i^N_m)}
g^{(i^r_l)}_{\La_r}(z)
\Big].
\nonumber
\end{eqnarray}
Then
\begin{eqnarray}
&&
G^{rl}(i^{N-1}_k)=
{\prod_{s=1}^{m-1}(i^{N-1}_k-i^N_s)
\over
i^{N-1}_k-i^N_m}
\Big[
1-
\prod_{s=1}^{m-1}{i^N_m-i^N_s \over i^{N-1}_k-i^N_s}
\prod_{s\neq l}^m{i^{N-1}_k-i^r_s \over i^N_m-i^r_s}
\Big].
\label{f42}
\end{eqnarray}
Using (\ref{f42}), (\ref{f23}) and (\ref{f22}) we have, from (\ref{f41}),
\begin{eqnarray}
&&
\df{\La}{i^N_m}-\df{i^r_l}{i^N_m}
\nonumber
\\
&&=
\sum_{k=1}^{m-1}
\prod_{j\in{\cal K},j\neq i^{N-1}_k}
{i^N_m-j \over i^{N-1}_k-j}
\Big[
1-
\prod_{s=1}^{m-1}{i^N_m-i^N_s \over i^{N-1}_k-i^N_s}
\prod_{s\neq l}^m{i^{N-1}_k-i^r_s \over i^N_m-i^r_s}
\Big]
\Big[
\df{\La}{i^{N-1}_k}-\df{i^{N-1}_{m}}{i^{N-1}_k}
\Big].
\label{f43}
\end{eqnarray}
Substituting  (\ref{f30}) into (\ref{f43}) we have (\ref{rel4}).

\end{document}